\documentclass[12pt,letterpaper]{amsart}
\usepackage{ amsmath,amsthm, amsbsy,amsfonts,amssymb, txfonts,tabularx,booktabs}
\usepackage[normalem]{ulem}
\usepackage{stmaryrd,mathrsfs,graphicx, amscd, tikz-cd, tikz, cancel}
\usepackage{tikz,enumitem}
\usetikzlibrary{matrix,arrows,decorations.pathmorphing}

\usepackage[active]{srcltx}
\allowdisplaybreaks
\numberwithin{equation}{section}

\usepackage[
	hypertexnames=false,
	hyperindex,
	pagebackref,
	pdftex,
	breaklinks=true,
	bookmarks=false,
	colorlinks,
	linkcolor=blue,
	citecolor=red,
	urlcolor=red,
]{hyperref}
\usepackage{hyperref}

\def\CC{{\mathbb C}}

\def\HH{{\mathbb H}}

\def\PP{{\mathbb P}}
\def\Pb{{\mathbb P}}
\def\QQ{{\mathbb Q}} 
\def\RR{{\mathbb R}} 
\def\Rb{{\mathbb R}}

\def\ZZ{{\mathbb Z}} 
\def\Zb{{\mathbb Z}}

\def\ssm{\smallsetminus}

\def\aff{{\rm aff}}

\def\G{\Gamma}
\def\g{\gamma}

\def\orb{{\rm orb}}

\def\prim{{\rm prim}} 
\def\int{{\rm int}} 
\def\el{{\rm ell}} 
\def\nod{{\rm nod}}

\def\Itwo{{\rm I_2}}

\def\bs{\backslash}
\newcommand{\eps}{\varepsilon}

\newcommand{\p}{\partial}

\def\Dcal{{\mathcal D}}

\def\Hcal{{\mathcal H}}

\def\Kcal{{\mathcal K}}
\def\Lcal{{\mathcal L}}
\def\Mcal{{\mathcal M}}
\def\Ocal{{\mathcal O}}
\def\Pcal{{\mathcal P}}

\def\Scal{{\mathfrak S}}
\def\Sscr{{\mathscr S}}
 
\def\Ucal{{\mathcal U}}

\def\Xcal{{\mathcal X}}

\def\Bscr{{\mathscr B}}  
\def\Cscr{{\mathscr C}}  
\def\Dscr{{\mathscr D}}

\def\la{\langle}
\def\ra{\rangle}

\def\pt{{\scriptscriptstyle\bullet}}

\newcommand\aut{\operatorname{Aut}}

\newcommand\Diff{\operatorname{Diff}}

\newcommand\Gr{\operatorname{Gr}}

\newcommand\jac{\operatorname{Jac}}

\newcommand\Mod{\operatorname{Mod}}

\newcommand\MW{\operatorname{MW}}

\newcommand\pic{\operatorname{Pic}}

\newcommand\SL{\operatorname{SL}}

\newcommand\supp{\operatorname{supp}}

\newcommand\Tcal{\operatorname{{\mathcal T}}}

\newcommand\Orth{\operatorname{O}}

\newcommand\PSL{\operatorname{PSL}}

\newcommand\Tr{\operatorname{Tr}}

\newtheorem{theorem}{Theorem}[section]
\newtheorem{lemma}[theorem]{Lemma}
\newtheorem{proposition}[theorem]{Proposition}
\newtheorem{corollary}[theorem]{Corollary}

\newtheorem{definition}{Definition}\numberwithin{definition}{section}

\theoremstyle{remark}

\newtheorem{remark}[theorem]{Remark}

\newtheorem{question}[theorem]{Question}

\title[Moduli of genus one fibered K3 surfaces]{Moduli spaces and period mappings of genus one fibered K3 surfaces}
\author{Benson Farb} 
\author{Eduard Looijenga}
\thanks{The first author was supported in part by National Science Foundation Grant No. DMS-181772 and the Eckhardt Faculty Fund.  The second author was  supported in  part by the Mittag-Leffler institute (Fall of 2021). Both authors are supported by the Jump Trading Mathlab Research Fund.}

\begin{document}

%\subjclass[2000]{81R10, 17B65}
%\keywords{}

\begin{abstract}
In this paper we construct various moduli spaces of K3 surfaces $M$ equipped with a surjective holomorphic map $\pi:M\to\Pb^1$ with generic fiber a complex torus (e.g., an elliptic fibration).  Examples include moduli spaces of such maps with primitive fibers; with reduced, irreducible fibers; equipped with a section.  Such spaces are closely related to the moduli space of Ricci-flat metrics on $M$.  We construct period mappings 
relating these moduli spaces to locally symmetric spaces, and use these to compute their orbifold fundamental groups.  

These results lie in contrast to, and exhibit different behavior than, the well-studied case of moduli spaces of polarized K3 surfaces, and are more useful for applications to the mapping class group $\Mod(M)$.    Indeed, we apply our results on moduli space to give two applications to the smooth mapping class group of $M$.

\end{abstract}
\maketitle

%\tableofcontents

\section{Introduction} 
A {\em K3 surface} is a closed, simply-connected complex surface $M$ admitting a nowhere vanishing holomorphic $2$-form.   Many K3 surfaces $M$ admit a {\em (holomorphic) genus one fibration}; that is, a surjective holomorphic map $\pi:M\to \Pb^1$ with finitely many singular fibers, and whose smooth fibers are Riemann surfaces of genus one.  Such fibrations play a central role in the theory of K3 surfaces.

Determining the image of the period mapping for a space $\mathcal M$ of complex structures on a compact manifold $M$, as well as the related 
problem of computing $\pi_1^{\rm orb}({\mathcal M})$, are fundamental but difficult problems (see, e.g.\ \cite{De}, \cite{K}  and \cite{CMP} for surveys/introductions).  In this paper we (mostly) solve these problems for 
genus one fibrations of K3 surfaces.  In summary: 

\begin{enumerate}
\item We construct the moduli spaces of various types of genus one fibered K3 surfaces, for example those with primitive fiber class; those with reduced, irreducible fibers; those equipped with a section; etc.  See Theorem \ref{theorem:existence}.
\item We construct period mappings relating these moduli spaces to locally symmetric varieties, and compute the images of these period mappings.  See Theorems \ref{thm:Mintell}, \ref{thm:nodal} and \ref{thm:Jac}.
\item We apply this to compute and relate the fundamental groups of these moduli spaces.  See Theorems \ref{thm:prim}, \ref{thm:Mintell}(3) and 
\ref{theorem:mint1}.
\item We apply the above to prove two results on the smooth mapping class group of $M$: a maps version of a theorem of Mosihezon (Theorem \ref{theorem:moishmaps}); and a Nielsen-type realization theorem (Theorem \ref{cor:MWNielsen}).
\end{enumerate}

\begin{remark}
In order for the moduli spaces we consider to be useful, particularly for applications to understanding mapping class groups of K3 surfaces, we want them to be Hausdorff.   Applications of this type also explain our focus on their orbifold fundamental groups.  For this reason we will exclusively deal with K3 surfaces endowed with a K\"ahler class. \emph{As a consequence,  most of our moduli spaces do not  have a complex structure}  and so this brings us in general outside the context of algebraic varieties (or stacks). This is in contrast with the well-studied  moduli spaces of polarized K3 surfaces, as these come with a natural structure of a locally symmetric variety (and  are more relevant for the study of symplectic mapping class groups).
\end{remark}

Associated to any genus one fibration of $M$ is the Poincar\'e dual $e\in H^2(M;\Zb)$ of the class of a generic fiber.  The vector $e$ is nonzero and, since it is the class of a fiber in a fiber bundle, is {\em isotropic}: $e\cdot e$=0.  The converse is also true: for any nonzero, primitive isotropic 
$e\in H^2(M)$ there is a complex structure on $M$ and a holomorphic genus one fibration with $e$ as the class of a generic fiber. In this way isotropic vectors $e\in H^2(M)$ will be central to this paper.

\subsection{Ricci-flat metrics and complex structures}\label{subsect:RF}
 Kodaira proved that all K3 surfaces are diffeomorphic.  In this paper $M$ denotes a closed 4-manifold in this diffeomorphism class; we give it the orientation for which the intersection pairing has signature $(3,19)$. According to Donaldson \cite{D}, $H^2(M; \Rb)$
comes with a natural \emph{spinor orientation}, that is,  an orientation of the bundle  of positive  subspaces of maximal dimension (so in this case that dimension is $3$) over the appropriate Grassmannian. 

In what follows the twistor construction plays a central role, so we recall some of the basics.  
Henceforth we assume that all metrics have unit volume. Given a Riemann metric $g$ on 
$M$, the space $P_g$ of its self-dual harmonic $2$-forms  defines a positive $3$-plane in $H^2(M; \Rb)$.  If $J$ is a complex structure on 
$M$ for which  $g$ is a K\"ahler metric, then the K\"ahler class $\kappa$ 
of this metric lies in $P_g$ and 
$\kappa\cdot\kappa =2$. The orthogonal complement of $\kappa$ in $P_g$ is a positive $2$-plane that inherits via the spinor orientation a canonical  orientation.  This determines a complex line $H^{2,0}(M,J)\subset H^2(M; \CC)$: if $(x,y)$ is an oriented orthonormal basis of the $2$-plane $P_g$ then $H^{2,0}(M,J)$ is spanned by $x+\sqrt{-1}y$; further, when we know 
$H^{2,0}(M,J)$ then we know the full Hodge decomposition of $H^{2}(M,J)$.  

In case  $g$ is Ricci flat there is a converse: each $\kappa\in P_g$ with $\kappa\cdot\kappa =2$ determines a complex structure $J$ on $M$ for which $g$ is a K\"ahler metric that has $\kappa$ as its K\"ahler class. This complex structure is unique and all complex structures $J$ for which $g$ is a K\"ahler metric so arise \cite{HKLR}; we then say that the pair $(g,J)$ is a \emph{K\"ahler-Einstein structure} on $M$. In particular, the set of K\"ahler-Einstein structures  with $g$ as underlying Ricci-flat metric is faithfully parametrized by the 2-sphere of radius $\sqrt{2}$ in $P_g$. In universal terms, if 
\[\Mcal_{\rm RF}:=\{\text{$\Diff(M)$-orbits of Ricci-flat metrics on $M$}\}\]
and if
\[ \Mcal_{\rm KE}:=\{\text{$\Diff(M)$-orbits of K\"ahler-Einstein structures on $M$}\},
\]
then once we have shown these to  be moduli spaces in our sense, the forgetful map 
$\pi:\Mcal_{\rm KE}\to \Mcal_{\rm RF}$ defined by $\pi(g,J):=g$ is an orbifold $S^2$-bundle: 
\begin{equation}
\label{eq:kemetrics1}
\begin{array}{rc}
S^2\longrightarrow& \Mcal_{\rm KE}\\
&\big\downarrow\pi\\
&\Mcal_{\rm RF}
\end{array}
\end{equation}

A starting point for this paper is the observation that a choice of nonzero isotropic vector $e\in H^2(M;\Zb)$ determines a smooth section of \eqref{eq:kemetrics1}, and that in fact 
the following stronger result is true.

\begin{proposition}[{\bf Metric $+$ isotropic vector $\rightsquigarrow$ complex structure}]
\label{proposition:KE2}
Let $M$ be the K3 manifold.  Fix a nonzero isotropic vector $e\in H^2(M;\Zb)$.  Then every 
Ricci-flat metric $g$ on $M$ determines a unique complex structure $J_g$ on $M$ for which $g$ is K\"ahler-Einstein and $e$ is the class of an 
effective divisor. This divisor can be chosen to be the sum of a smooth genus $1$ curve and a nonnegative linear combination of smooth genus $0$ curves.  
\end{proposition}

See Section \ref{proof:prop:KE2} for a proof.  Additional conditions are needed in order that $e$ be represented by a 
smooth genus $1$ curve. But if that is the case, then $e$ defines a basepoint-free linear system of dimension $1$ (a copy of $\PP^1$), i.e., it gives rise to a genus one fibration $\pi: (M,J)\to S$ 
(where $S$ is a copy of $\PP^1$) such that $e$ is the Poincar\'e dual of the class of a fiber and (hence) $e$  is the image of the orientation class of $S$ under $\pi^*: H^2(S)\to H^2(M)$.

\subsection{Existence theorems}  Proposition \ref{proposition:KE2} allows us to form a number of moduli spaces of holomorphic genus one fibrations as complex structures and fiber classes vary. The first of these is: 
\begin{equation}\label{eqn:Mprim}
\Mcal_\prim:=
\left\{
\begin{array}{l}
\Diff(M)\text{-orbits of triples $(g, J, \pi)$, where $(g,J)$ makes $M$ a}\\
\text{K\"ahler-Einstein K3 surface for which  $\pi$ is a holomorphic }\\
\text{genus one fibration of $M$ with primitive fiber class}
\end{array}\right\}
\end{equation}

A singular fiber in a genus one fibration of $M$ is {\em integral} if it is reduced and irreducible; equivalently, it is nodal (Kodaira type $\rm{I}_1$) or cuspidal (Kodaira type $\rm{II}$).   
We will call a  genus one fibered structure on $M$ {\em integral}  if each singular fiber is integral.  This is a property that in the situation of Proposition \ref{proposition:KE2} can easily be characterized homologically: there should not exist an $\alpha\in H^2(M)$ with $\alpha\cdot\alpha=-2$ that is perpendicular to both $e$ and to $\kappa$.
The second moduli space we shall consider is : 
\[
\begin{array}{ll}
\Mcal_\int:=&\text{the locus in $\Mcal_\prim$ for which the fibration $\pi$ is integral}
\end{array}
\]

Not every holomorphic genus one  fibration on  $M$ admits a holomorphic section.  When it does, it is called an {\em elliptic fibration}, since that section chooses a basepoint for each fiber, under which each smooth fiber becomes an elliptic curve with this basepoint as the identity group element (each reduced singular fiber also becomes a group, depending on the type of fiber--see below).  Two sections of a genus one fibration differ by a holomorphic fiberwise translation, although  that isomorphism will in general not 
preserve a given K\"ahler class. But in case the holomorphic genus one fibration is integral and a 
section with class $\sigma$ is given, then  there is a canonical choice for the K\"ahler class, namely  
$\sqrt{1/2}(e+3\sigma)$ ($e+3\sigma$ is  the class of  an ample line bundle of degree 4, see Corollary \ref{cor:polarization}), and hence the class of a (unique) K\"ahler-Einstein metric. This allows us to regard 
\begin{equation}\label{eqn:diagram}
\begin{array}{ll}
\Mcal_\int^{\el}&:=\{\Diff(M)\text{-orbits of integral elliptic fibrations on $M$}\}
\end{array}
\end{equation}
as a subset  of $\Mcal_\int$. 
There is also a map going the other way: 
\[
\jac : \Mcal_\int\to \Mcal_\int^\el.
\]
It assigns to a genus one  fibration its {\em Jacobian fibration}, replacing fiberwise each cubic curve by its Jacobian (see \S\ref{section:jacobian} for more details) and makes $\Mcal_\int^\el$ a  retract of  $\Mcal_\int$.  The  Jacobian fibration is defined for a general genus one fibration of $M$,
but since that Jacobian fibration  does not seem to have a natural K\"ahler class that varies smoothly in families (as in the integral case), we refrain from introducing  a moduli space $\Mcal^\el_\prim$.

The discriminant of a genus one fibration $M\to S$ with only integral fibers yields a divisor $D$ on $S$ of degree $24$, with a nodal fiber contributing with multiplicity one  and a cuspidal fiber with multiplicity two. It is the same as for its Jacobian fibration. Let us write $\Mcal^{\le 2}_{(24)}$ for the $\PSL_2(\CC)$ orbit space of positive degree 24 divisors on $\PP^1$ with multiplicities $\le 2$ and $\Mcal_{(24)}\subset \Mcal^{\le 2}_{(24)}$  for the locus of reduced divisors. Let $\Mcal_\nod\subset \Mcal_\int$ be the set of genus one fibrations all of whose singular fibers are nodal. 

The following theorem summarizes the moduli spaces considered in this paper and 
various natural maps between them.  As indicated below, some of these spaces (e.g., $\Mcal_{\rm RF}$) have previously been considered in the literature.  Among our main results will be to give uniformizations of most of these moduli spaces.
 
\begin{theorem}[{\bf Existence theorem}]
\label{theorem:existence}
With the notation above, the diagram 
\begin{center} 
\label{eq:relatingpi1}
\begin{tikzcd}\label{eq:relatingpi1}
\Mcal_\nod\arrow[d, "\jac"] \arrow [r, hook ]  & \Mcal_\int\arrow[d, "\jac"] \arrow [r, hook]  & \Mcal_\prim \arrow [r]\arrow[rd]  & \Mcal_{\rm KE}\arrow [d ]\\
\Mcal^\el_\nod\arrow [u,  bend left=15, hook ]\arrow [d]\arrow [r, hook ]& \Mcal_\int^\el \arrow [d]\arrow [u,  bend left=15, hook ]& &\Mcal_{\rm RF}\\ 
\Mcal_{(24)} \arrow [r, hook ]& \Mcal^{\le 2}_{(24)}
\end{tikzcd}
\end{center}
is one of moduli spaces (and morphisms between them) in the sense of the conventions on moduli spaces in \S\ref{conv:moduli} below.  Each of these moduli spaces is connected.
\end{theorem}

We will abuse notation and often identify each of these moduli spaces with its base.  
With the exception of $\Mcal_\nod$ and $\Mcal_\nod^\el $ we will give fairly 
concrete descriptions of these moduli spaces, often in terms of arithmetic groups acting on open subsets of 
homogeneous spaces.  This will allow us to discuss the orbifold fundamental groups of these moduli spaces.  
We will see that, apart of $\Mcal_{(24)}$ and $\Mcal^{\le 2}_{(24)}$ where it is evident, only $\Mcal_\int^\el$ 
(and hence $\Mcal_\nod^\el$) live in the quasi-projective category.  That moduli space is called in \cite{loenne2021}  the \emph{Miranda moduli space}.   

\begin{remark}
Before we proceed, a word on terminology. We want to use the notions of `orbifold', `orbifold fundamental group' and `moduli space' in such a manner that we are for instance able to say that the moduli space of elliptic curves can be given as the orbifold $\SL_2(\ZZ)\bs\HH$ and that therefore its orbifold fundamental group is $\SL_2(\ZZ)$ (which in this case also happens to be the mapping class group of a torus). The appropriate  language is that of Deligne-Mumford stacks in the smooth category, but for us the more elementary conventions stated in \S\ref{conv:moduli} already do the job.
\end{remark}

\subsection{Topology of the moduli spaces} 
The following table gives a quick, incomplete summary of the main results of this paper.

\renewcommand{\arraystretch}{1.2}
\begin{table}[h]
\centering
\caption{A summary of some of the results of this paper.}
\begin{tabularx}{\linewidth}{@{} *4{>{\raggedright\arraybackslash}X} @{}}
\toprule
  Moduli space
  & Uniformization
  & $\pi^\orb_1$
  & Reference\\
  
\midrule
$\Mcal_{\rm RF}$ & $\G\bs\Gr^+_3(H_\RR)^\circ$ & $\G$ &\S\ref{subsection:Teich} and Thm.\ref{thm:prim}\\
$\Mcal_{\rm KE}$ & $\G\bs E(H_\RR)$& $\G$ & \S\ref{subsection:Teich} and Thm.\ref{thm:prim}\\
$\Mcal_\prim$ & see \S\ref{section:prim2}& $\G_e$& Thm.\ref{thm:unipotent radical}\\
$\Mcal_\int^\el$ & $\G(e)\bs\Gr^+_2(H_\RR)^\circ$&$\widetilde\G(e)$&  Thm.\ref{thm:Mintell}-(3)\\
$\Mcal_\int$ &$T^{20}\times [0,\infty)$-bundle over $\G(e)\bs\Gr^+_2(H_\RR)^\circ$& $\G_e\times_{\G(e)}\widetilde\G(e)$ &Thm.\ref{theorem:mint1}\\  
  \bottomrule
\end{tabularx}
\end{table}

The orbifold fundamental groups of the moduli spaces above are related to a certain arithmetic groups, as we now explain.  
Let $H$ denote $H^2(M;\Zb)$ equipped with its intersection form.  This is an even unimodular lattice of signature $(3,19)$ and these properties characterize $H$ up to (isometric) isomorphism.  The intersection form extends to $H_\Rb:=H\otimes\Rb$.  The orthogonal group $\Orth(H)$ is an arithmetic lattice in $\Orth(H_\Rb)$, 
a Lie group isomorphic to the real semisimple Lie group $\Orth(3,19)$.  

Let $\Mod(M):=\pi_0(\Diff(M))$ be the smooth mapping class group of $M$.  The action of $\Diff(M)$ on $M$ induces a representation 
\[\rho:\Mod(M)\to\Orth(H)\]  
whose image, which we will denote by $\G$, is {\it a priori} contained in the index $2$ subgroup $\Orth^+(H)\subset \Orth(H)$ preserving the spinor orientation, but is in fact known to be equal to that group \cite{B}. 

Let $e\in H$ be a primitive isotropic vector. 
 As we explain in more detail in \S\ref{subsection:lattices1} below, the lattice 
\[
H(e):=e^\perp/\Zb e
\]  
is even, unimodular and of signature $(2,18)$, and these properties characterize its isomorphism type. The natural action of the $\G$-stabilizer $\Gamma_e$ of $e$ on $H(e)$ induces a representation 
\[
\Gamma_e\to \Orth(H(e))
\]
whose image has index $2$ (it is defined by a spinor orientation), and will be denoted by $\Gamma(e)$.  This gives a (noncanonically  split) short exact sequence
\[
0\to H(e)\xrightarrow{E} \G_e\to \G(e)\to 1
\]
where $H(e)$ can be identified with the unipotent radical of $\Gamma_e$, consisting of those elements of $\Gamma_e$ acting trivially on $H(e)$. 

\begin{theorem}[\textbf{Topology of $\Mcal_\prim, \Mcal_{KE},$ and $\Mcal_{RF}$}]\label{thm:prim}
The map on orbifold fundamental groups induced by the diagram 
\[\Mcal_\prim\to \Mcal_{KE}\to \Mcal_{RF}\] 
is (up to conjugacy) naturally isomorphic to 
\[\Gamma_e\subset \G =\G.\] 
\end{theorem}

The new part of Theorem \ref{thm:prim} concerns $\Mcal_\prim$; the other claims are recalled in \S\ref{subsection:Teich}.

In what follows, we choose a fixed elliptic fibration $\pi: M\to S$ with only nodal fibers as singular fibers. The discriminant $D$ of $\pi$  is then a $24$ element subset of $S$. We use this fibration as a basepoint for the various moduli spaces introduced here. 

The group of isotopy classes of orientation-preserving diffeomorphisms  of the pair $(S,D)$ is the  spherical braid group.  That group is generated by the set of elementary braids. An \emph{elementary braid} is given by an arc in $S$ that connects two distinct points of $D$ but whose interior avoids $D$  (it is rather an isotopy class of such)
and  the associated  elementary braid is a half Dehn twist whose support is contained in a regular neighborhood of that arc, so that it exchanges the two points of $D$.  It is not hard to show that an orientation-preserving diffeomorphism of $(S,D)$ lifts to $M$. But in order that it represents a loop in a moduli space of genus one fibrations, it must lift in a fiber preserving manner, and this is not the case unless special conditions (that pertain to the restriction of $\pi: M\to S$ to that arc) are met. 

The symmetric space of the orthogonal group $\Orth(H(e)_\RR)$ is the Grassmannian of positive $2$-planes in $H(e)_\RR$, denoted $\Gr_2^+(H(e)_\RR)$. Let 
\[\Gr^+_2(H(e)_\RR)^\circ\subset \Gr^+_2(H(e)_\RR)\] denote the locus of positive $2$-planes that have no $(-2)$-vector (that is, a vector with self-intersection $-2$) in their orthogonal complement. It is of course $\G (e )$-invariant.

The following theorem gives a uniformization via a period mapping of the moduli space $\Mcal_\int^\el$, uses this to compute $\pi_1^{\rm orb}(\Mcal_\int^\el)$, and gives a modular interpretation of certain degenerations in $\Mcal_\int^\el$.   The first item is likely 
known to experts.

\begin{theorem}[{\bf Topology of $\Mcal_\int^\el$}]\label{thm:Mintell}
The orbit space $\G(e)\bs\Gr^+_2(H(e )_\RR)$  is a locally symmetric quasi-projective orbifold and 
\begin{enumerate}
\item there is an explicit period mapping (a morphism of orbifolds)
\[
P: \Mcal_\int^\el \to \G(e)\bs\Gr^+_2(H(e )_\RR)
\]
that is an open embedding with image $\Gamma(e)\backslash\Gr^+_2(H(e)_\RR)^\circ$, the complement of an irreducible, locally symmetric hypersurface. 
\item The general point of this hypersurface represents  a genus one fibration of $M$ with a fiber of  Kodaira type $\Itwo$ (\footnote{Recall that such a fiber is the union of two copies of $\Pb^1$ intersecting transversally in $2$ points.}). 

\item The period map $P$ exhibits $\pi_1^{\orb}(\Mcal_\int^\el)$ as a group $\widetilde\G(e)$ in an extension of groups
\begin{equation}
\label{eq:exact11}
 1\to \pi_1(\Gr^+_2(H(e )_\RR)^\circ)\to  \widetilde\G(e)\to  \Gamma(e)
\to 1,
\end{equation}
where $\pi_1(\Gr^+_2(H(e )_\RR)^\circ)$ is normally generated  by the square of an elementary spherical  braid in $S$ defined by an arc connecting two discriminant points along which 
there is  a degeneration into a singular fiber of Kodaira type ${\rm I}_2$-fiber. 
\item The natural fiberwise-involution  of the universal  elliptic 
fibering, which can be regarded as an element of the inertia 
subgroup of $\pi_1^{\orb}(\Mcal_\int^\el)$, maps to minus the identity in $\Gamma(e)$. 
 
\end{enumerate}
\end{theorem}

The following theorems is a computation of the oribifold fundamental group of $\Mcal_\int$. 

\begin{theorem}[{\bf Topology of $\Mcal_\int$}]
\label{theorem:mint1}
There is a diagram of   split short exact sequences
\begin{center} \label{eq:relatingpi1}
\begin{tikzcd}
\label{eq:relatingpi1}
0 \arrow [r] &H(e)\arrow [r]\arrow [d, equal] &\pi^\orb_1(\Mcal_\int)\arrow [r]\arrow [d] & \pi^\orb_1(\Mcal_\int^\el)\arrow [r]\arrow [d] & 1\\
0\arrow [r] &H(e)\arrow [r, "E"] & \Gamma_e\arrow [r] & \G(e)\arrow[r] & 1
\end{tikzcd}
\end{center}
In particular, the right hand square is cartesian, so that 
\[\pi^\orb_1(\Mcal_\int)\cong \G_e\times_{\G(e)}\widetilde\G(e).\]
\end{theorem}

We will prove the following theorem in \S\ref{subsect:proof of ref{thm:nodal}}.

\begin{theorem}[{\bf Topology of $\Mcal_\nod$}]\label{thm:nodal}
The inclusion  $\Mcal_\nod\subset \Mcal_\int$  induces a  surjection on orbifold fundamental groups  whose kernel is normally generated by the  lift to $M$ of a third power  of an elementary braid  defined by an arc in $S$ connecting  two nodal fibers along which there is 
a degeneration into a cuspidal fiber.
\end{theorem}

\begin{remark}\label{rem:}
We have not been able to obtain for $\Mcal_\nod$ a concrete description of the same type as we have for the other moduli spaces. This is because we do not know how read off from the Hodge structure of a genus one fibration on $M$ whether it has a cuspidal fiber.
\end{remark}

\subsection{Universal Jacobian fibrations} Our next main result is a modular interpretation of Theorem \ref{theorem:mint1}, which we now explain.  Associated to any K3 surface $X$ endowed with 
holomorphic integral genus $1$ fibration $\pi:X\to S$ is a {\em Jacobian fibration}  
\[
\jac(\pi) \to S,
\]  
which is an elliptically fibered K3 surface that has the same base and the same discriminant as $\pi$, but replaces  each smooth fiber with its Jacobian. The fibration $\jac(\pi)\to S$ comes with a holomorphic section,  namely its zero section, whereas $\pi$ need not have one; indeed, the two fibrations are (fiberwise) isomorphic if and only if $\pi$ admits a section.   This construction globalizes to a `universal Jacobian'  map 
\[\jac: \Mcal_{\rm int}\to \Mcal_\int^\el\]
converting holomorphic genus $1$ fibrations to elliptic fibrations.  See \S\ref{section:jacobian} for more details. 

\begin{theorem}[{\bf The universal Jacobian construction}]\label{thm:Jac}
In the diagram of moduli spaces
\begin{center}
\begin{tikzcd}
\Mcal_\int\arrow [r, "\jac"] & \Mcal_\int^\el\cong \Gr^+_2(H(e)_\RR)^\circ\arrow[l, bend right=30, hook]
\end{tikzcd}
\end{center}
the  forgetful map appears as a section of $\jac$.  
The map $\jac$ has the structure of a $C^\infty$-fiber bundle of orbifolds whose  fibers have the homotopy type of an $20$-dimensional torus. More precisely, the torus $T(e):=\RR/\ZZ\otimes H(e)$ endowed with its $\G(e)$-action defines 
a torus bundle over $\Mcal_\int^\el$  (namely the orbit space of 
$T(e)\times \Gr^+_2(H(e)_\RR)^\circ$ with respect to the diagonal action of $\G (e)$), 
and there is a fiberwise homotopy equivalence (as orbifolds) over $\Mcal_\int^\el$  of $\Mcal_\int$ with this torus bundle.
In particular, the induced maps on $\pi_1^{\rm orb}$ induce the long exact sequence of Theorem \ref{theorem:mint1}, as well as its splitting.
\end{theorem}

\subsection{Two applications to smooth topology} The above results on moduli spaces of holomorphic genus one fibrations have applications to the smooth mapping class group $\Mod(M):=\pi_0(\Diff(M))$.  If we are given a particular structure $\Sscr$ on $M$ 
(such as an elliptic fibration  with only nodal fibers as singular fibers), then
there is an associated mapping class group $\Mod_\Sscr(M)$, defined as the connected component group of the group of diffeomorphisms of $M$ that preserve this structure. 
$\Mod_\Sscr(M)$ comes with a forgetful homomorphism  
$\Mod_\Sscr(M)\to \Mod(M)$.  The connection with the moduli space $\Mcal_\Sscr$ comes from the fact that if $\Mcal_\Sscr$ is connected then
 the monodromy of universal bundle over $\Mcal_\Sscr$ induces a representation 
\[\pi_1^{\rm orb}(\Mcal_\Sscr)\to \Mod_\Sscr(M).\]
The period mappings constructed in this paper produce homomorphisms $\pi_1^{\rm orb}(\Mcal)\to \Gamma$ 
that factor through $\Mod(M)$ and whose image we can often determine.  For example, the connectedness of $\Mcal_{\rm nod}$ implies the following (see \S\ref{section:applics} for details).  

\begin{theorem}[{\bf Moishezon for maps}]
\label{theorem:moishmaps}
Given any $\gamma\in O(H)^+$ fixing a nonzero, isotropic vector $e\in H_2(M;\Zb)$ there exists:
\begin{enumerate}
\item A complex structure $J$ on $M$;
\item A holomorphic (with respect to $J$) elliptic fibration $\pi:M\to\Pb^1$ whose fibers have 
homology class $e$.
\item A fiber-preserving (with respect to $\pi$) diffeomorphism $f:M\to M$ such that 
$f_*=\gamma$.
\end{enumerate}
\end{theorem}

Items 1 and 2 of Theorem \ref{theorem:moishmaps} are due to Moishezon (see (\cite{FM}, Cor.\ 7.5)), and give a new proof of that result for K3 surfaces.

As a second application of results we will deduce the following (see Theorem \ref{thm:nielsenMW} below for a more precise statement).

\begin{corollary}[Nielsen realization for $H(e)$]\label{cor:MWNielsen}
Given a genus one fibration of a $M$ with only nodal fibers and fiber class $e$, then $H(e)$, when regarded as a subgroup of the $\G_e$, lifts to a group of fiber-preserving diffeomorphisms of $M$.
\end{corollary}

Corollary \ref{cor:MWNielsen} explains  how the abelian group $R_u(\G_e )\cong H(e )$ of rank 20 appears in algebraic geometry as a monodromy group: if we fix an elliptic K3 surface $\pi: X\to\PP^1$ with section $\sigma$, then it is the monodromy group of the family of 
genus one fibered K3 surfaces whose Jacobian fibration is isomorphic to the pair $(\pi, \sigma)$.  
The group of diffeomorphisms given in Corollary \ref{cor:MWNielsen} can be thought of 
as a Mordell-Weil group of rank $20$ in the smooth category, where we note  that  the maximal  rank in the holomorphic category is at most $18$.  A subsequent paper will be devoted to these  mapping class groups.

\subsection{Conventions on orbifold groups and orbifold structures on moduli spaces}
\label{conv:moduli}
In this paper, an {\em orbifold}  appears always as a global quotient, i.e., as an orbit space $\G\bs T$ of a smooth manifold $T$ by a group $\G$ acting properly discontinuously by diffeomorphisms on $T$.  In case $T$ is simply-connected,  we declare the {\em orbifold fundamental group} of $\G\bs T$ to be  $\pi_1^{\rm orb}(\G\bs T):=\Gamma$. Note that we here allow the action to be non-faithful; we call its kernel, necessarily finite,  the \emph{inertia subgroup} of $\pi_1^{\rm orb}(\G\bs T)$. If $T$ is only connected and $\tilde T\to T$ is a universal covering with Galois group $\pi_1(T)$, then the set of all lifts of all elements of $\G$ to $\tilde T$ is 
a group $\tilde \G$ that is  an extension $\tilde\G$  of $\G$ by  $\pi_1(T)$. This group acts properly discontinuously on $\tilde T$ and we regard the evident bijection  
$\tilde\G\bs \tilde T\to \G\bs T$ as an isomorphism of orbifolds, so that $\tilde\G$  is the orbifold fundamental group of $T$.

We use a similar convention for moduli spaces. These will exclusively concern a class of structures (denoted $\mathscr{S}$) that we can put on a manifold $M$, where we assume $\Sscr$ closed under pullback by a diffeomorphisms of $M$. Then {\it grosso modo} a moduli space for $\Sscr$-structures on $M$ puts an orbifold structure on the set of 
$\Diff(M)$-equivalence classes in $\Sscr$. To be precise, consider orbifolds $\G\bs T$  with $T$ simply connected  that parametrize  diffeomorphism classes of $\Sscr$-structures on  $M$ as follows: we are given  a smooth fiber bundle $\Ucal\to T$ with fiber diffeomorphic to 
$M$ and with an $\Sscr$-structure given on each fiber in a smoothly varying manner. We also assume that the $\G$-action on $T$ has been lifted to  $\Ucal$ in a way that preserves this structure. 
Then we say  that such a family is {\em universal} up to $\G$-action if for every family $\Ucal'\to T'$ of this 
type (so with $T'$ simply-connected, but here no group $\G'$ acting on it is assumed) fits in cartesian diagram 
\[
\begin{CD}
 \Ucal'@>>>\Ucal\\
 @VVV  @VVV\\
T' @>>> T\\
\end{CD}
\]
with the pair  of horizontal maps  being structure preserving  and being unique for this property up 
to postcomposition with an element of $\G$. In particular (take $T'$ a singleton) a structure-preserving isomorphism between two fibers of $\Ucal\to T$ is then always induced by an element of $\G$. It is  not hard to show that $\Ucal\to T$ is unique  
up to the $\G$-action (one may think of this as an almost final object of a category: it is unique up to $\G$-action). We  then say that $\G\bs(\Ucal\to T)$  (or simply $\G\bs T$, when the remaining data are understood)  is the \emph{moduli space} for manifolds diffeomorphic with 
$M$ and endowed with an $\Sscr$-structure.

\subsection{Connection with spherical braids}
\label{subsection:discrmap}

We conclude this introduction with a final remark and question.  The closure of the image of   $\Mcal_\int^\el\to \Mcal^{\le 2}_{(24)}$ is the set of degree 24 divisors (so with multiplicities $\le 2$)  that can be given by the sum of a cube and a square, i.e, that lie in  a linear system generated by $3D_0$ and $2D_1$ with $D_0$  and $D_1$  effective divisors of degree $8$ resp.\ $12$ (see for instance the discussion in $\S$ 2.1 of \cite{HL}). 
If we denote that locus by $\Dscr^{\le 2}_{(24)}\subset \Mcal^{\le 2}_{(24)}$, then Hasimoto-Ueda \cite{HU} prove that the resulting  map  $\Mcal_\int^\el\to \Dscr^{\le 2}_{(24)}$ is of degree one. This map is also easily seen to be quasi-finite (i.e., has finite fibers).  

L\"onne determined in \cite{loenne2021} the induced map on 
 orbifold fundamental groups $\pi_1^\orb(\Mcal_\int^\el)\to \pi_1^\orb(\Mcal^{\le 2}_{(24)})$  
 (where the latter is also known as the spherical braid group).  The map  $\Mcal_\int^\el\to \Dscr^{\le 2}_{(24)}$ is not proper, because we do not allow Kodaira fibers of type $\mathrm{I}_2$, 
 which also have discriminant multiplicity $2$. We can however show  that the restriction to reduced discriminants:
\[
\Mcal_\nod^\el\to \Dscr_{(24)}=\Dscr^{\le 2}_{(24)}\cap \Mcal_{(24)}
\] 
gives a local homeomorphism that is also a normalization. Yet we do not know the answer to 
the following.

\begin{question}
What is $\pi_1^\orb(\Dscr_{(24)})$? What is its image in $\pi_1^\orb(\Mcal_{(24)})$?
\end{question}

\subsection*{Acknowledgements}
It is a pleasure to thank Dan Margalit, Curt McMullen and Yunpeng Zi for comments and corrections on an earlier draft of this paper. We also thank the anonymous referees for their extensive comments and corrections.

\section{Preliminary material}
In this section we present some preliminary material that will be used throughout the paper.

\subsection{Period mappings for moduli spaces of Ricci-flat metrics}
\label{subsection:Teich}
We now recall the construction of $\Mcal_{\rm RF}$ and $\Mcal_{\rm KE}$, as well as their associated period mappings and use this to determine their fundamental groups.

Fix a K3 manifold $M$.  As recalled in the introduction, $M$ is oriented  for which $H:= H^2(M;\Zb)$, endowed with the intersection pairing,  is an even unimodular lattice of signature $(3,19)$  and $H_\RR= H^2(M;\Rb)$ comes with a spinor orientation. Both are preserved  by the action of $\Diff(M)$ on $H$ and its image  $\G$  is  all of  $\Orth(H)\cap \Orth^+(H_\RR)$. This is subgroup of index two in $\Orth(H)$ that does not contain minus the identity, and so has trivial center. In particular, $\G$ acts faithfully on the Grassmannian $\Gr^+_3(H_\RR)$ of positive $3$-planes in $H_\RR$, 
which we recall is the symmetric space of $\Orth(H_\RR)$. 

\begin{definition}[{\bf Torelli space $\Tcal(M)$ of Ricci-flat metrics}]
Let $\Tcal(M)$ be the space of isometry 
classes of K3 manifolds $X$ equipped with a  Ricci-flat metric $g$ of unit volume and an {\em $H$-marking}, that is, an isomorphism $H^2(X)\to H$ which preserves the intersection pairing and the spinor orientation.  The space $\Tcal(M)$ is called the {\em Torelli space} of $M$.
\end{definition}

Let $\Diff_H(M)$ denote the kernel of the representation $\Diff (M)\to \Orth(H)$, 
so that  $\Diff(M)/\Diff_H(M)$ can be identified with $\Gamma$.  The space 
$\Tcal(M)$  can then be characterized as the  $\Diff_H(M)$-orbit space of the space of unit volume, Ricci-flat metrics on $M$. It comes with a natural $\G$-action.
%It is a  simply-connected $57$-dimensional manifold on which the  group 
%$\Gamma$ acts properly discontinuously.

\begin{remark}
Another option would be to work with the Teichm\"uller space, defined  as a $\Diff(M)^0$ orbit space of the space of unit volume, Ricci-flat metrics on $M$.   Each connected component of that space is isomorphic to $\Tcal(M)$ and  the {\em Torelli group} $\Diff_H(M)/\Diff(M)^0$ permutes them  simply transitively. It  is still an open question whether the Torelli group of $M$ is trivial.
 \end{remark}

Let $\Gr^+_3(H_\RR)$ denote the Grasmannian of positive $3$-planes in $H_\RR$.  
Assigning to a metric $g$  on $M$ the space $P_g$ of its self-dual harmonic $2$-forms on $(M,g)$
defines a $\G$-equivariant {\em period mapping} 
\begin{equation}
\label{eq:pmapping1} %P:
 \tilde{P}:\Tcal(M)\to \Gr^+_3(H_\RR). 
\end{equation}
The Torelli Theorem for K\"ahler-Einstein K3 surfaces (see for instance \cite{looijenga:hyp})  asserts that $\tilde{P}$ is an open embedding with image 
\[\Gr^+_3(H_\RR)^\circ\subset \Gr^+_3(H_\RR)\] 
the set of positive $3$-planes in $H_\RR$ that have no $(-2)$-vector of $H$ in their orthogonal complement.  The set of such planes is a $\Gamma$-invariant, locally finite union of codimension $3$ submanifolds, so that $\Gr^+_3(H_\RR)^\circ$ is open in $\Gr^+_3(H_\RR)$.
The induced map
\[P:\Mcal_{\rm RF}\to \Gamma\bs\Gr^+_3(H_\RR)^\circ\]
is a diffeomorphism of orbifolds. Now, $\Gr^+_3(H_\RR)$ is the symmetric space corresponding to the semisimple Lie group ${\rm O}(3,19)(\Rb)$, and so it is nonpositively curved and contractible.  Since $\Gr^+_3(H_\RR)^\circ$ is the complement of codimension $3$ submanifolds, it follows that $\pi_1(\Gr^+_3(H_\RR)^\circ)=0$, and so there is an isomorphism
\[P_*:\pi_1(\Mcal_{\rm RF})\stackrel{\cong}{\to}\Gamma,\]
where the left-hand side is an orbifold fundamental group.

It is shown in \cite{looijenga:hyp} that $\Tcal(M)$ supports a family 
\[\Ucal_{\Tcal(M)}\to \Tcal(M)\cong \Gr^+_3(H_\RR)^\circ\] 
of K3-manifolds endowed with a  unit volume, Ricci-flat metric to which the $\G$-action lifts. This family has the universal property of the conventions on moduli spaces in \S\ref{conv:moduli}.  

This in turn leads to a universal $H$-marked family of K3 surfaces endowed with a unit volume, Ricci-flat K\"ahler  (or K\"ahler-Einstein) metric as follows.  The twistor construction shows that for a given Ricci-flat metric on $M$, the set of the complex structures on $M$  for which this metric is a K\"ahler metric is a 2-sphere, and if we do this universally we find a 
$2$-sphere bundle $p_E: E(M)\to \Tcal(M)$ such that the pull-back 
\[p_E^*\Ucal_{\Tcal(M)}\to E(M)\]
 yields the universal family in question. 
 
The space $E(M)$ can also be described in terms of the period map: if  a Ricci-flat metric on $M$ defines the  positive $3$-plane $P\subset H_\RR$, then the $2$-sphere can be identified with the sphere of radius $\sqrt{2}$ in $P$. Indeed, the imaginary part of the K\"ahler metric is a closed $2$-form whose class is a $\kappa_P\in P$ with self-product $2$ (this reflects the fact that the metric is unital). We denote the corresponding  2-sphere bundle by $E(H_\RR)\to \Gr^+_3(H_\RR)$, so that $E(M)$ gets identified with $E(H_\RR)^\circ$, the restriction of  $E(H_\RR)$ to $\Gr^+_3(H_\RR)^\circ$. 

This entire picture is $\G$-equivariant, and so the map descends to a diffeomorphism
\[\Mcal_{\rm KE}\to \Gamma\bs E(M),\]
making $\Mcal_{\rm KE}$ into a moduli space in the sense of the conventions on moduli spaces in \S\ref{conv:moduli}.  Since $p_E$ is an $S^2$-bundle, the
projection  $\G\bs E(H_\RR)\to \G\bs \Gr^+_3(H_\RR)^\circ$
induces an isomorphism on orbifold fundamental groups, so that
 \[(p_E)_\ast:\pi_1(\Mcal_{\rm KE}) \to \pi_1(\Mcal_{\rm RF})\cong\Gamma\]
 is an isomorphism.

\subsection{Groups and lattices attached to a primitive isotropic vector}
\label{subsection:lattices1}
This subsection is purely group-theoretic.  It gives the structure of the stabilizer in the arithmetic group $\Gamma$ of an isotropic vector.  This will be crucial for describing the fundamental groups of various moduli spaces of holomorphic genus $1$ fibrations of $M$.

\emph{In the rest of this paper we fix a primitive isotropic vector} $e \in H$, so with $e\cdot e=0$.   It does not matter which one we choose because all such vectors belong to the same $\G$-orbit. This follows from the following well-known fact about lattices 
(see for instance \cite{LP}): if $L$ is an even lattice, which such as $H$,  has a copy  $\mathbf{U}\perp\mathbf{U}$ as a direct summand (here $\mathbf{U}$ stands for the lattice $\ZZ^2$ endowed with the quadratic form $(x,y)\in \ZZ^2\mapsto xy$), then its orthogonal group acts transitively on the primitive vectors of a given length. Any such vector is represented in one of the $\mathbf{U}$-summands and then it follows (by exchanging the basis vectors of the other copy of $\mathbf{U}$) that  the subgroup $\Orth^+(L)$ that preserve each spinor orientation will have the same property. 

We  next make some observations regarding  the $\G$-stabilizer $\G_e $.  This stabilizer    leaves invariant the short flag 
$\{0\}\subset \ZZ e \subset e ^\perp\subset H$. We shall write $H_e $ for $e ^\perp$ and $H(e )$ for $e ^\perp/\ZZ e$. We may identify the latter with the image of $H_e$ under the map $e\wedge : H\to \wedge^2 H$ so that $H(e)\cong e\wedge H_e$.
The real vector space $H(e )_\RR$ has signature $(2,18)$ and inherits from $H$ a spinor orientation: 
its bundle of  positive  defintite $2$-planes comes with an orientation. 

Since $H/H_e $ can be identified with the dual  of $\ZZ e $, the group $\G_e $ acts trivially on it and so the unipotent radical  $R_u(\G_e )$ of $\G_e $ consists of the elements of $\G_e $ that act also trivially on $H(e )$.  We denote the image of $\G_e$ in the orthogonal group of $H(e )$ by $\G(e )$.

The following proposition collects some useful properties of these groups and lattices.

\begin{proposition}[{\bf Properties of $\Gamma_e$ and related lattices}]
\label{proposition:elattice}
With the notation as above, the following hold:
\begin{enumerate}
\item The lattice $H(e )$ is even  unimodular of signature $(2,18)$ and comes with a natural spinor orientation.  
\item  $\G(e )$ is the index $2$ subgroup of $\Orth (H(e ))$  that preserves this spinor orientation.  
\item  There is 
a (noncanonically split) exact sequence
\[
0\to H(e )\xrightarrow{E} \G_e \to \G(e )\to 1
\]
that identifies $H(e )$ with the unipotent radical $R_u(\G_e )$ of $\G_e $. Here we represent $\g\in H(e )$ as a $2$-vector $e \wedge\g\in \wedge^2 H$ with $\gamma\in H_e$ as above and $E$ assigns to the latter the {\em Eichler transformation} 
\begin{equation}\label{eqn:eichler}
E(e \wedge \g): x\in H_\QQ \mapsto x + (x\cdot e )\g- (x\cdot \g)e -\tfrac{1}{2}(\g\cdot \g)(x\cdot e )e \in H_\QQ
\end{equation}
\item  $ \G(e )$ contains the central element $-1$ and acts transitively on the set of all primitive vectors in $H(e )$ of a given self-product. 
\item The $(-2)$-vectors in $H_e $ make up a $\G_e $-orbit, and the group of Eichler transformations that we identified with $H(e )$ acts transitively on the set of $(-2)$-vectors in $H(e )$ that lie over a given $(-2)$-vector of $H(e )$.
\end{enumerate}
\end{proposition}
 
\begin{proof}
Choose $e '\in H$ such that $e '\cdot e =1$. Upon replacing $e '$ by $e '-\tfrac{1}{2}(e'.e')e$ we can and will assume that $e '\cdot e '=0$ so that $(e ,e ')$ spans a copy $U$ of $\mathbf{U}$.
Its orthogonal complement $U^\perp$  is then an even  unimodular lattice of signature $(2,18)$ and hence isomorphic to $\mathbf{E_8}(-1)\perp \mathbf{E_8}(-1)\perp \mathbf{U}\perp \mathbf{U}$. It is clear that $U^\perp$ maps isomorphically onto $H(e )$.
It is known that the reflections in the $(-2)$-vectors of such a lattice generate the  index $2$ subgroup  of its orthogonal group  that preserves spinor orientations and that the  primitive vectors of a given self-product make up a single orbit under that group.  So $H(e )$ has that property. 
This also implies that $ \G(e )=\Orth^+(H(e ))$. That the latter contains $-1$ is clear.

We next exhibit the exact sequence. Any orthogonal transformation of $H_\QQ$ which fixes $e $ and acts trivially on $H(e )_\QQ$ is of  the form
\[
E(e \wedge \g): x\in H_\QQ \mapsto x + (x\cdot e )\g- (x\cdot \g)e -\tfrac{1}{2}(\g\cdot \g)(x\cdot e )e \in H_\QQ
\]
for some $\g\in H_{e,\QQ}$ (this indeed only depends on $e \wedge \g$, that is, only depends on the image $\g'$ of $\g$ in $H(e)_\QQ$). Note that  if $x\in H_{e,\QQ}$, then
$E(e \wedge \g)$ takes $x$ to $x- (x\cdot \hat\g)e $. So if we ask that this transformation preserves the lattice $H$, then we must have that
$x\cdot \g\in \ZZ$ for all $x\in H_e $. Since $H(e )$ is unimodular, this amounts to $\g'\in H(e )$.  But this necessary condition clearly also suffices. We thus identify the kernel $\G_e \to \G(e )$ with $H(e )$. The subgroup $\G_e \cap\G_{e '}\subset \G_e $  maps isomorphically onto $\G(e )$ and so provides a splitting of the displayed exact sequence.

The $(-2)$-vectors in $H(e )$ form a single $\G(e )$-orbit,  and so it remains to prove that if $\alpha,\alpha'\in H_e $  are $(-2)$-vectors with the same image in $H(e )$, then there exists an Eichler transformation of the above type that takes $\alpha$ to $\alpha'$. Clearly $\alpha-\alpha'=n e $ for some $n\in\ZZ$.  Since $H_e$ is unimodular  there exists a $\g\in H_e $ is such that  $\alpha\cdot \g=1$, and then   $E(e \wedge n\g)$ takes $\alpha$ to $\alpha-ne =\alpha'$.
\end{proof}

\section{Moduli spaces of elliptic K3-surfaces: the proof of Theorem \ref{thm:Mintell}}

In this section we prove Theorem \ref{thm:Mintell}.

\subsection{The integrality criterion} 

We begin with some basic facts that we will need; we refer to \cite{LP} for the assertions stated here without proof.   

Let $X$ be a K3 surface. The space $H^{1,1}(X)$, which is the orthogonal complement of $H^{2,0}(X)\oplus\overline{H^{2,0}(X)}$ in $H^2(X;\CC)$, is defined over $\RR$, and the intersection form restricted to $H^{1,1}(X; \RR)$ has signature $(1,19)$. 
The spinor orientation on $X$ singles out a connected component  $\Cscr_X^+$ of the space of $x\in H^{1,1}(X; \RR)$  with $x\cdot x>0$: it is the connected component that contains all the K\"ahler classes.  Let $\kappa$ be one such class.

A class in $H^2(X;\Zb)$ is representable by a divisor if and only if it lies in $H^{1,1}(X)$. The linear equivalence class of that divisor is then unique since $H^1(X, \Ocal_X)=0$, so that
\[
\pic (X)= H^2(X;\Zb)\cap H^{1,1}(X)= H^2(X;\Zb)\cap H^{2,0}(X)^\perp.
\] 
If $v\in \pic (X)$ is nonzero and such that $v\cdot v\ge 0$ or $v\cdot v=-2$, then either $v$ or $-v$ is the class of an effective divisor, depending on the sign of $v\cdot\kappa$. The effective divisors thus obtained generate the whole semigroup $\pic^+(X)\subset \pic (X)$ of effective divisor classes.

A special role is played by the classes of smooth curves on $X$ of genus $0$ and $1$, often called \emph{nodal} resp.\ \emph{elliptic} classes. The adjunction formula shows that the self-intersection number is then $-2$ resp.\ $0$. The semigroup $\pic^+(X)$ is already generated by the nodal classes, the primitive elliptic classes and $\pic(X)\cap \Cscr_X^+$.  

Each nodal class $\alpha$ defines an orthogonal  reflection $s_\alpha$ in $H^2(X;\Zb)$ defined by
$x\mapsto x+(\alpha\cdot x)\alpha$. It preserves both $\pic(X)$ and $\Cscr^+_X$.
The set of all reflections in nodal classes generate a Coxeter subgroup  $W_X\subset \Orth(H)$ (which together with this generating set make it a Coxeter system) that acts as such 
on $H^{1,1}(X; \RR)$ with `fundamental chamber' the cone 
\[
C_X:=\{x\in H^{1,1}(X; \RR)\, :\, x\cdot \alpha\ge 0 \text{ for every nodal class $\alpha$}\}.
\]
The $W_X$-orbit of $C_X$ contains $\Cscr_X^+$ and so $C_X\cap\Cscr_X^+$ is a fundamental domain for the action of $W_X$ on $\Cscr_X^+$.  The following is well known (\cite{LP}):

\begin{lemma}[{\bf Roots in $\pic(X)$}]
\label{lemma:kcone}
Every $\alpha \in \pic (X)$  with $\alpha\cdot \alpha=-2$ is a \emph{root} for $W_X$: it lies in the $W_X$-orbit of a nodal class and $\alpha$ can be written as  a $\ZZ$-linear combination of nodal classes for which all nonzero coefficients have the same sign 
(and we call accordingly $\alpha$ a \emph{positive} or a  \emph{negative} root; is also the sign of the function $x\mapsto x\cdot \alpha$  on the interior of $C_X$).

A class in $H^2(X; \RR)$  is a  K\"ahler class if and only it lies in  $\Cscr_X^+$ and has positive intersection product with every nodal class (this is equivalent to it lying in the interior of $C_X\cap\Cscr_X^+$). Any class in $\Cscr_X^+$ not fixed by a reflection in a root lies in the $W_X$-orbit of a K\"ahler class. 
\end{lemma}

A  primitive elliptic class $e$ defines a dominant weight: it lies in $C_X$. It also
defines a $1$-dimensional linear system $|e|$, and this linear system is a genus one fibration of $X$ over a copy of $\PP^1$. The irreducible components of the reducible fibers of this fibration  are nodal and have zero intersection number with $e$. 
This has an important implication, namely that if $d$ is the class of an 
effective divisor, then $d\cdot e\ge 0$ with the equality  sign implying that $d$ is a nonnegative combination of irreducible components of fibers of $|e|$. 

The discussion above implies the following. 

\begin{proposition}[{\bf Integrality criterion}]
\label{prop:isotropic class}
Let  $e\in \pic^+(X)$ be primitive with $e\cdot e=0$. Then the following hold:
\begin{enumerate}
\item $e$ is representable as a sum of an elliptic class $e'$ plus a nonnegative linear combination of nodal classes, and lies in the $W_X$-orbit of $e'$.
\item $e$ is the class of a genus one fibration if  and only if $e$ is dominant weight, that is  $e\cdot\alpha\ge 0$ for every nodal class.
\item The fibration has only integral fibers if and only if $e$ is strictly dominant in the sense $e\cdot\alpha>0$ for every nodal class. 
\end{enumerate}
\end{proposition}

We also have the following.
\begin{corollary}\label{cor:polarization}
Let $e$ be the class of a genus one fibration with only integral fibers. If $\sigma$ is the class of a section, then $\kappa:=3e+\sigma$ is the class of a K\"ahler metric.
\end{corollary}
\begin{proof}
We verify that $\kappa$ satisfies the conditions of Lemma \ref{lemma:kcone}, i.e.,  that 
$\kappa\in \Cscr^+_X$ and that $\kappa$ has positive intersection product with every nodal class.

We first observe that $\sigma$ is represented by a smooth genus zero curve and hence nodal.  
Thus $\sigma\cdot \sigma=-2$, and so $\kappa\cdot \sigma=3\cdot 1-2=1>0$.  Since   $\kappa\cdot e=1$, it the follows that $\kappa\cdot \kappa=3\cdot 1+1=4$. So $\kappa\in \Cscr^+_X$.
If $\alpha$ is a nodal class different from $\sigma$, then 
$\alpha\cdot e > 0$ and $\alpha\cdot \sigma\ge 0$ and hence  $\kappa\cdot\alpha > 0$. 
\end{proof}

\begin{remark}\label{rem:}
The class $\kappa$ is in fact the class of an ample divisor. It defines on $X$ a basepoint-free linear system $|\kappa|$ of dimension $3$, but is not very ample: the resulting map $X$ to $\PP^3$ has as image a twisted cubic $Y$ (a smooth rational curve of degree $3$) such that the map $X\to Y$  is the given  elliptic fibration. However, $3\kappa$ is very ample.

\end{remark}

\subsection{K\"ahler-Einstein metrics: proof of Proposition \ref{proposition:KE2}}

Fix a primitive isotropic vector $e\in H$.  Suppose $M$ is endowed with a Ricci-flat metric  $g$ with associated  positive oriented $3$-plane $P\subset H_\RR$.  As mentioned above, every $\kappa\in P$ with $\kappa\cdot\kappa=2$ determines a complex structure $J$ on $M$ for which the given metric is K\"ahler and has $\kappa$ as its class. The Hodge structure of this  complex structure can be recovered from  the \emph{oriented}
2-plane $\Pi :=\kappa^\perp\cap P$ by means of the recipe described in Subsection \ref{subsect:RF}.

The spinor orientation on $H(e )_\RR$ implies that  each  positive  $2$-plane $\Pi\subset H(e )_\RR$ defines a Hodge structure $H_\Pi^{\pt,\pt}$ on $H(e )$ that is polarized by the given pairing: if $(x,y)$ is an oriented  orthonormal basis of $\Pi$, then $H^{2,0}_\Pi$   is the $\CC$-span of
$x+\sqrt{-1}y$, $H^{0,2}_\Pi$ its complex conjugate and $H^{1,1}_\Pi$ is the complexification of $\Pi^\perp$. This gives the  symmetric space  $\Gr_2^+(H(e )_\RR)$ the structure of a bounded symmetric domain. Since $\G(e )$ is a finite index subgroup of the orthogonal group of $H(e )$, the Baily-Borel theory tells us  that the  orbit space $\G(e )\bs \Gr_2^+(H(e )_\RR)$ is then in a natural way a quasi-projective variety.

Similarly, a  positive  $2$-plane $\tilde\Pi$ in $H_{e,\RR}$ defines a Hodge structure (but no longer polarized) on $H_e$. The projection $H_e\to H(e )$ becomes a morphism of Hodge structures if we endow 
$H(e )$ with the Hodge structure defined by the image of $\tilde\Pi$. 

\begin{proof}[Proof of Proposition \ref{proposition:KE2}]
\label{proof:prop:KE2}
There is precisely one $\kappa\in P$ with $\kappa\cdot\kappa=2$ such that the linear form
$x\in P\mapsto \kappa\cdot x$ is a positive multiple of  $x\in P\mapsto  e\cdot x$.  This $\kappa$  lifts $g$ to a K\"ahler-Einstein structure $(g,J)$ with the property that $e$ is perpendicular to $H^{2,0}(M,J)$. It follows that $e$ is of type $(1,1)$ and $\kappa\cdot e>0$. Proposition \ref{prop:isotropic class} then tells us that $e$ is  the class of an effective 
 divisor of the type stated here.
\end{proof}

\subsection{Periods of elliptic fibrations: proof of Theorem \ref{thm:Mintell}}
With all of the above in hand, we can now prove Theorem \ref{thm:Mintell}.

\medskip
\noindent
{\bf Proof of Item (1). }
Let us first make explicit the period mapping in question.  Fix a $\sigma\in H$ with $e\cdot \sigma=1$ and $\sigma\cdot\sigma=-2$. The classes $\{e,\sigma\}$ span a copy $U\subset H$ of  the basic hyperbolic lattice $\mathbf{U}$.  Further, $U^\perp$, which is contained in $e^\perp=H_e$, maps isomorphically onto $H(e)$. 

Suppose that $X$ is an elliptically fibered 
$K3$-surface with only integral fibers.  This gives a fiber class $e_X$ and a section class $\sigma_X$ spanning a copy 
$U_X\subset H^2(X)$ of $\mathbf{U}$ as above. Then there exists a marking 
$H^2(X;\Zb)\stackrel{\cong}{\to} H$ that 
takes $e_X$ to $e$ and $\sigma_X$ to $\sigma$. This marking is induced by a diffeomorphism and we use that  diffeomorphism to pull back all structure on $X$ to $M$, so that $M$ gets a complex structure 
$J$ for which $e$ is the class of an elliptic, integral  fibration that has $\sigma$ as the class of its section. The orbit of this structure under the action of the stabilizer of $\{e,\sigma\}$ in $\Diff(M)$ on $M$ independent of our choices. 

The Hodge structure of $(M,J)$ is given by an oriented, positive 2-plane in $U^\perp_\RR$, which therefore  is given by  an oriented positive 2-plane $\Pi\subset U^\perp_\RR$. By Proposition
\ref{prop:isotropic class}, $\Pi^\perp \cap H_e$ does not contain any root. Hence if $\Pi'\subset
H(e)_\RR$ is the image of $\Pi$, then $\Pi'{}^\perp \cap H(e)$ does not contain a $(-2)$-vector.

The same argument works for families of such surfaces with a simply-connected base $T$: we then obtain a holomorphic map $T\to\Gr^+_2(H(e)_\RR)^\circ$ that is unique up to postcomposition with an element of $\G(e)$.  There is thus an induced 
period map 
\[P: \Mcal_\int^\el \to \G(e)\bs\Gr^+_2(H(e )_\RR)^\circ.\] 

The Torelli theorem implies that $P$ is an isomorphism of quasi-projective varieties.

Let $\alpha\in H(e)$  be such that $\alpha\cdot\alpha=-2$. 
Then 
\[
\Gr_2^+(\alpha^\perp_\RR)=\Gr_2^+(H(e)_\RR)^{s_\alpha}
\]
 parametrizes the set of those Hodge structures in $H$ for which $\alpha$ is of type $(1,1)$. A generic point of $\Gr_2^+(\alpha^\perp_\RR)$ will have the property that the $(1,1)$-part is the spanned $\{e, \sigma, \alpha\}$. It is representable by a genus one fibration $\pi: M\to S$ for which $e$ is the fiber class, $\sigma$ the class of a  section and $\alpha$ is the class of a  divisor (a root). 
 
By part (4) of Proposition \ref{proposition:elattice}, the group $\G(e)$  acts transitively on the set  
$\alpha\in H(e)$ with $\alpha\cdot\alpha=-2$. So the deleted locus in $\G (e)\bs\Gr_2^+(H(e)_\RR)$ 
(the hypersurface mentioned in the theorem) is irreducible. 

\medskip
\noindent
{\bf Proof of Item (2). }

We must show that a general element of this hypersurface represents  a genus 
one fibration of $M$ with a fiber of  Kodaira type $\Itwo$.  First note that if the elliptic fibration on $M$ 
(with fiber class $e$ and section class $\sigma$) comes with a type $\Itwo$ fiber, then the irreducible components of this fiber are nodal curves.  The section intersects only one of these curves, and so if $\alpha$ is the class of 
the other, then $\alpha$ is a nodal class with $\alpha\cdot \sigma=\alpha\cdot e=0$. Identifying $H(e)$ with 
the orthogonal complement of $\ZZ e+\ZZ\sigma$ in $H$ (so that $\alpha$ is now regarded as an element of 
$H(e)$), then it is clear that this elliptic fibration defines an element of $\Gr_2^+(\alpha^\perp_\RR)$. 

The Torelli theorem shows that the converse also holds: given a $(-2)$-class $\alpha\in H$ perpendicular to both $e$ and $\sigma$, then there exists a complex structure $J$ on $M$ for which the Picard group $\pic(M,J)$ is spanned by $\{e, \sigma, \alpha\}$, where $e$ is the class of an elliptic fibration with section class $\sigma$.  Replacing $\alpha$ by $-\alpha$ if necessary, 
we can assume that $\alpha$ is effective. An effective divisor representing $\alpha$ then lies in a finite union of fibers. Since the section represented by $\sigma$ has zero intersection with this divisor,  the fibers in question are all reducible. The class of an irreducible component is perpendicular to both $e$ and $\sigma$ and hence a multiple of $\alpha$. It follows that $\alpha$ is represented by a nodal curve $C$ contained in a fiber $F$.  

The same argument shows that $F$ has no other irreducible components beside $C$ and the irreducible component $C'$ that meets the section. Since the class of $C'$ is $e-\alpha$, we have $C'\cdot C'=-2$ and so $C'$ is a nodal curve. It also follows that $C\cdot C'=2$. By means of  small deformation of the complex structure we can arrange that $C$ and $C'$ meet transversally so that we  get a fiber of Kodaira type $\Itwo$. 

\medskip
\noindent
{\bf Proof of Item (3). }
The identification of orbifold fundamental groups
$\Mcal_{\rm ell}^{\rm int}\cong\G(e)\bs\Gr^+_2(H(e)_\RR)^\circ$ is now clear. In particular, the kernel of 
\[
\pi^\orb_1(\G(e)\bs\Gr^+_2(H(e)_\RR)^\circ)\to \pi^\orb_1(\G(e)\bs\Gr^+_2(H(e)_\RR))\cong \G(e)
\]
is normally generated by the boundary in $\Gr^+_2(H(e)_\RR)^\circ$ of a small closed disk  
$\varphi: D^2\hookrightarrow\Gr^+_2(H(e)_\RR)$ that is $s_\alpha$-invariant and which meets $\Gr^+_2(\alpha^\perp_\RR)$ transversally with $\varphi^{-1}\Gr^+_2(\alpha^\perp_\RR)=\{0\}$.
Then we find over $D^2$  a degenerating family of elliptic fibrations with central fiber 
an elliptic fibration with just one non-nodal fiber, that fiber being of type  $\Itwo$. This 
descends to a map from the $s_\alpha$-orbit space of  $D^2$ to $\G(e)\bs\Gr^+_2(H(e)_\RR)$ that is transversal 
to the image of $\Gr_2^+(\alpha^\perp_\RR)$ so that the restriction to the boundary represents a simple loop around this hypersurface.  The resulting map $\la s_\alpha\ra\bs\p D^2\to \Scal_{24}\bs\Mcal_{0,24}$ 
yields a simple braid of the asserted type, with $\p D^2$ naturally lifting to $\Mcal_{0,24}$.

\medskip
\noindent
{\bf Proof of Item (4). } It remains to see that the fiberwise involution $\iota$ of an elliptic fibration $M\to \PP^1$ acts on $e^\perp/e$ as minus the identity. Or equivalently, that $\QQ e+\QQ\sigma$ is the fixed-point  set of $\iota$ in $H^2(M;\QQ)$. For this it is helpful to recall that an integral elliptic fibration is locally over its base in Weierstra\ss\ form: 
it is given as 
\[y^2= x^3+a(s)x+b(s)\] with $s$ a local coordinate on the base.  In these coordinates the fiberwise involution is given by 
\[\iota(x,y,s)=(x,-y,s).\] This shows that the  orbit space $M_\iota$ of  
the elliptic fibration has the structure of a $\PP^1$-bundle over $\PP^1$ (the local chart becomes $(x,s)$) that comes with a section. The space  $H^2(M_\iota; \QQ)$ is spanned by the class of the section and a fiber 
and these are sent under the natural map $H^2(M_\iota; \QQ)\to H^2(M; \QQ)$  to $2\sigma$ and $2e$ respectively.
As the image of this map is $H^2(M; \QQ)^\iota$, the assertion follows.

\begin{remark}\label{rem:}
The quotient  $\Mcal_\int^\el:=\G(e )\bs\Gr^+_2(H(e )_\RR)^\circ$ 
is the moduli space  that has been considered  by L\"onne in \cite{loenne2021} (for the value $d=4$ in that paper), the  main result being a presentation of its orbifold fundamental group. 
We here find a description that brings it in relation to $\G (e )$ in a way that is similar
 to how a braid group relates to a symmetric group: over each $(-2)$-reflection in 
 $\G(e )$ lies an element of infinite order represented by a simple braid in $\G (e )\bs \Gr^+_2(H(e )_\RR)^\circ$.
\end{remark}

\section{Topology of $\Mcal_\prim$: proof of Theorem \ref{thm:prim}}
\label{section:prim2}

In this section we prove Theorem \ref{thm:prim}.

\subsection{$\Mcal_\prim$ and its topology}
Recall  that we have fixed a primitive isotropic vector $e\in H$.  The locus in $\Tcal (M)$ for which $e$ is the class of a genus one fibration (resp.\ an integral genus one fibration) will be denoted by $\Tcal(M)_e$ (resp.\ $\Tcal(M)_e^\circ$). 

We have seen that the period map gives $\G$-equivariant diffeomorphism $\Tcal (M)\cong \Gr_3^+(H_\RR)^\circ$. In terms of this isomorphism we can characterize $\Tcal(M)_e$ and $\Tcal(M)_e^\circ$ as subloci of  $\Gr_3^+(H_\RR)^\circ$. Proposition  \ref{proposition:KE2}
associates to every Ricci-flat metric $g$ on $M$ a complex structure $J$ for which $g$ is a K\"ahler metric (whose class we denote by $\kappa$) and $e$ is the class of an effective divisor of a particular type. In order that  $e$ be the class of a genus one fibration resp.\  an integral genus one fibration a necessary and sufficient condition is that for every positive root $\alpha$ we have $\alpha\cdot e\ge 0$ resp.\ 
$\alpha\cdot e> 0$. These properties can be entirely stated in terms of 
the positive $3$-plane $P\subset H_\RR$ defined by $g$, and so this gives the corresponding 
characterizations.

The following theorem tells  us more about $\Tcal(M)_e$. Except for the claim that  it is 
simply-connected,  its assertions  can be found, at least implicitly, in the literature.  It implies Theorem \ref{thm:prim}.  We will give its proof in \S\ref{proof:unip rad} below.

\begin{theorem}\label{thm:unipotent radical}
The locus $\Tcal(M)_e \subset  \Tcal(M)$ for which $e$ is the class of a holomorphic 
genus one fibration of $M$ is open 
and  $\G_e $-invariant. The  restriction of the universal bundle $\Ucal_{\Tcal(M)}\to\Tcal(M)$ to 
$\Tcal(M)_e $ comes with a $\G_e $-invariant 
factorization
\begin{center}
\begin{tikzcd}[column sep=scriptsize]
\Ucal_{\Tcal(M)_e }\arrow[dr] \arrow[rr]{} & &\Pcal_{\Tcal(M)_e} \arrow[dl]\\
&\Tcal(M)_e
\end{tikzcd}
\end{center}
where $\Pcal_{\Tcal(M)_e }\to \Tcal(M)_e $ is a $\PP^1$-bundle, thus exhibiting 
the universal property of this restriction: it is the universal family of $H$-marked K\"ahler-Einstein K3 surfaces  endowed with a genus one fibration with class $e$.
In particular,  the period map defines an isomorphism of orbifolds 

\[\Mcal_\prim\cong \Gamma_e\bs \Tcal(M)_e.\] 
Furthermore, $\Tcal(M)_e$ is simply connected, so that there is an isomorphism 
\[\pi_1(\Mcal_\prim)\cong \Gamma_e.\]
\end{theorem}

\subsection{Proof of Theorem \ref{thm:unipotent radical}}
\label{subsection:uniprad}

The construction underlying Proposition \ref{proposition:KE2} suggests that we consider the map 
\[
\begin{array}{rl}
h':\Gr_3^+(H_\RR)&\to \Gr_2^+(H_{e,\RR})\\
P&\mapsto P\cap H_{e, \RR}
\end{array}
\]
and the map
\[\begin{array}{rl}
h'':\Gr_2^+(H_{e,\RR})&\to\Gr_2^+(H(e)_\RR)\\
\Pi&\mapsto \text{image of $\Pi$ in $H(e)_\RR$}.
\end{array}
\]

Note that the  positive-definite $2$-planes in $H_{e,\RR}\subset H_\RR$ and $H(e)$ are  naturally oriented, so that they define complex lines in the complexifications of $H$, $H_e$ and $H(e)$.  That line is the $(2,0)$-part of a  weight two Hodge structure  on its ambient  vector space. The maps $h'$ and $h''$ can be understood in this context  as passing to a Hodge substructure and to a Hodge quotient structure. The link with Proposition \ref{proposition:KE2} is that if $[P]\in \Gr_3^+(H_\RR)$ is defined by a Ricci-flat metric on $M$, then $h'([P])$ gives the Hodge structure of the associated K\"ahler-Einstein structure.  Observe that $h'$ and $h''$ are $\G_e$-equivariant (where $\G_e$ acts on $\Gr_2^+(H(e)_\RR)$ through $\G (e)$). 

\begin{lemma}\label{lemma:factorization1}
The projection $h'':\Gr_2^+(H_{e,\RR})\to\Gr_2^+(H(e)_\RR)$ lives in the holomorphic category: the target is a Hermitian symmetric domain and $h''$ has the structure of  a complex affine line bundle over that domain.
\end{lemma} 
\begin{proof}
It is convenient to identify $\Gr_2^+(H_{e,\RR})$ with an open subset of the quadric in 
$\PP(H_{e,\CC})$ defined by the pairing; the point associated to $\Pi$ being the image of the isotropic vector
$z=x+\sqrt{-1}y$, where $(x,y)$ is an orthonormal oriented basis of 
$\Pi$. Then the fiber of $h''$ passing through 
$[z]\in \PP(H_{e,\CC})$ is the complex affine line parametrized by $\lambda\in\CC\mapsto [z+\lambda e]$. 
%It is easy to check that $H(e)_\RR$ acts transitively on such a fiber. 
\end{proof}

\begin{lemma}\label{lemma:factorization2}
The map  $h':\Gr_3^+(H_\RR)\to \Gr_2^+(H_{e,\RR})$
is a fiber bundle of real-hyperbolic spaces of dimension $19$.
\end{lemma} 
\begin{proof}
Let $\Pi$ be  a positive-definite $2$-plane  in $H_{e,\RR}$.  We determine the preimage of $[\Pi]\in \Gr_2^+(H_{e,\RR})$ under $h'$. The orthogonal complement of $\Pi$ in $H_\RR$, $\Pi^\perp$, has signature $(1,19)$ and  $h'{}^{-1}[\Pi]$ is identified with the set of $\kappa\in \Pi^\perp$ with $\kappa\cdot \kappa=2$ and $\kappa\cdot e>0$. This is indeed a real-hyperbolic space (the symmetric space for the orthogonal group of $\Pi^\perp$). 
\end{proof}

We will be interested in the restriction of the composite 
\[
h:=h''\circ h': \Gr_3^+(H_\RR)\to\Gr_2^+(H(e)_\RR)
\] 
to $\Tcal (M)_e\subset \Tcal (M)\cong\Gr_3^+(H_\RR)^\circ$.
By Lemma  \ref{lemma:factorization1},  $h''$ has contractible fibers and contractible domain. This is why our  focus will be on $h'$, or rather, its restriction to $\Tcal (M)_e$.

We first choose a positive, oriented 2-plane $\Pi$ in $H_{e,\RR}$ and  determine the  preimage of $[\Pi]\in \Gr_2^+(H_{e,\RR})$  in $\Tcal(M)$,  $\Tcal(M)_e$  and $\Tcal(M)^\circ_e$; these preimages will appear as subsets of the real-hyperbolic space given by Lemma  \ref{lemma:factorization2}, i.e., the set $\Hcal_\Pi$ of $\kappa\in \Pi^\perp$ with $\kappa\cdot\kappa=2$ and $\kappa\cdot e>0$. 
For this we shall invoke a bit of the Coxeter-Vinberg theory of reflection groups acting in such a hyperbolic space.

We denote  the set of $(-2)$-vectors in $H\cap \Pi^\perp$ by $R_\Pi$ and  refer  to its elements as 
\emph{$\Pi$-roots}. For any $\Pi$-root $\alpha$, the associated orthogonal reflection $s_\alpha$  in $H$  
leaves $\Pi$ pointwise fixed and $\Hcal_\Pi$ invariant. These reflections generate a Coxeter subgroup 
$W(R_\Pi)$ of $\G$.  The reflection hyperplanes in $\Hcal_\Pi$ make up  locally finite collection of subsets
and  decompose $\Hcal_\Pi$ into (relatively open) \emph{faces}. Such a  face is a  convex, locally polyhedral cone (it may not be definable by a finite number of inequalities and so it need not be a polyhedral cone).
An open face is  called a \emph{$W(R_\Pi)$-chamber}; this is also a connected component of the complement 
$\Hcal_\Pi^\circ$ of the union of the reflection hyperplanes. The group $W(R_\Pi)$ permutes the  chambers simply transitively. 

\begin{lemma}\label{lemma:} Let $\kappa\in\Hcal_\Pi$.
A necessary and sufficient condition that the $3$-plane $\Pi+\RR\kappa$ is an element of $\Gr_3^+(H_\RR)^\circ\cong \Tcal(M)$ is that  $\kappa\in\Hcal_\Pi^\circ$. So this identifies $h'{}^{-1}[\Pi]\cap \Tcal(M)$ with  $\Hcal_\Pi^\circ$.
\end{lemma}
\begin{proof}
It is clear that $\kappa\in\Hcal_\Pi^\circ$ is equivalent to the orthogonal complement of $\Pi+\RR\kappa$ not containing a $(-2)$-vector, i.e., to $\Pi+\RR\kappa$ giving an element of 
$\Gr_3^+(H_\RR)^\circ$. So any such $\kappa$ is the class of a  K\"ahler-Einstein structure for which $e$ represents an effective divisor.
\end{proof}

So if  $\kappa\in \Hcal_\Pi^\circ$, then there exists a K\"ahler-Einstein structure on $M$ for which $\kappa$ is the K\"ahler class and for which the oriented 2-plane $\Pi$ defines the Hodge structure. For any $\alpha\in R_\Pi$, either  $\alpha\cdot \kappa>0$ or $\alpha\cdot \kappa<0$; in the first case $\alpha$ is the  class of an  effective divisor. The simple $\Pi$-roots $\alpha$ relative to the chamber containing $\kappa$  (i.e., for which $\alpha\cdot\kappa>0$, but for which $\alpha$ cannot be written as the sum of two $\Pi$-roots with that property) are represented by smooth irreducible $(-2)$-curves on $M$. Furthermore $e$ is effective for that complex structure; it is in fact the class of an elliptic fibration plus a nonnegative linear combination of $(-2)$-curves.
We need to impose further restrictions on $\kappa$ in order that $e$ is the class of an elliptic fibration.

Let us denote the $\Pi$-roots $\alpha$ that are  perpendicular to $e$ by $R_{\Pi,e}$. The reflection in
the roots $R_{\Pi,e}$ generate a Coxeter  subgroup $W(R_{\Pi,e})\subset W(R_\Pi)$  which transitively permutes the chambers that have the improper point $[e]$ defined by $e$ in their closure.  Let $\Kcal_\Pi$ denote the union of the faces in $\Hcal_\Pi$  that have $[e]$ in their closure. This subset of $\Hcal_\Pi$ can be defined 
by the following (possibly infinitely many) strict inequalities: for all $\Pi$-roots $\alpha$ with $e\cdot\alpha > 0$, we have  $\kappa\cdot\alpha > 0$. We can restrict ourselves to the collection $\Phi_e$  of all the $\Pi$-roots $\alpha\in R_\Pi\ssm R_{\Pi,e}$ that define a $W(R_\Pi)$-chamber having $[e]$ as an improper point.  So this collection is invariant under the  $W(R_\Pi)$-stabilizer  $W(R_\Pi)_e$ of $e$. Indeed,  if we fix a  $W(R_\Pi)$-chamber which has $e$ as an improper point, then $\Phi_e$ is the  $W(R_\Pi)_e$-orbit of the simple roots relative to this chamber that are not perpendicular to $e$. It is  clear that $\Kcal_\Pi$  defines  a convex cone $\Hcal_\Pi$ (the intersection of $\Hcal_\Pi$ with a convex cone). %The theory of Coxeter and Vinberg tells us that $\Kcal_\Pi$ is an open in  $\Hcal_\Pi$ and is in fact a chamber 
%for the Coxeter subgroup $W^e(R_\Pi)\subset W(R_\Pi)$ generated by the reflections in the roots in $\Phi_e$.

The union $\Kcal_\Pi^\circ$ of open faces (=$W(R_\Pi)$-chambers) contained in $\Kcal_\Pi$ consists of the $\kappa\in \Hcal_\Pi$  with the property that  
for a root $\alpha$, the inequality  $e\cdot\alpha > 0$ implies   $\kappa\cdot\alpha >0$.  This translates into:  if  $\kappa\in\Hcal_\Pi$, then  $\Pi+\RR\kappa$ is defined by a genus one fibration with $\kappa$ as its K\"ahler class  if and only if $\kappa\in  \Kcal_\Pi^\circ$. In other words, this identifies 
$h'{}^{-1}[\Pi]\cap \Tcal(M)_e$ with  $\Kcal^\circ_\Pi$.  If $\kappa\in\Kcal_\Pi^\circ$ and $\alpha\in R_\Pi$ is a simple root relative to $\kappa$ so that  $\alpha$ is the class of a $(-2)$-curve, then this curve lies in a fiber of the genus one fibration if and only if $\alpha\cdot e=0$. Otherwise this curve  is class of a multi-section.
In particular, this genus one fibration has only integral fibers if and only if there are no roots $\alpha$ with $\alpha\cdot e=0$. In other words, this is the case if and only if $\Kcal_\Pi^\circ=\Kcal_\Pi$ (so $\Phi_e$ is then a root basis of $R_\Pi$ and consists of  the  classes of multi-sections of the fibration.)
\\

Recall that $\Tcal(M)\cong \Gr^+_3(H_\RR)^\circ$  is obtained from $ \Gr^+_3(H_\RR)$ by removing the fixed point sets of the reflections in $(-2)$-vectors in $H$. These are all of codimension $3$ and that is why $\Tcal(M)$ is simply-connected. In the proof of Corollary \ref{cor:sc} below, we  shall use a similar argument for $\Tcal(M)_e$.

\begin{corollary}\label{cor:sc}
The locus $\Tcal(M)_e$ is simply-connected.
\end{corollary}
\begin{proof}

The role of  $\Gr^+_3(H_\RR)$ is played here by a  subset $U_e\subset\Gr^+_3(H_\RR)$ parametrizing  positive $3$-planes $P$ with the property that for every $(-2)$-vector $\alpha \in H$ with $\alpha\cdot e>0$, we also have that $\alpha\cdot e_P>0$, where  $e_P$ is the orthogonal projection of $e$ in $P$. 
The local finiteness of the collection of codimension $3$-loci $\Gr^+_3(H_\RR)^{s_\alpha}$ implies that $U_e$ is open. We claim  that $U_e$ is also contractible. To see this,  note that the projection $U_e\to \Gr^+_3(H_\RR)$
factors as
\begin{center}
\begin{tikzcd}
U_e \arrow[rr, hook]\arrow[dr]& &\Gr^+_3(H_\RR)\arrow[dl]\\
 &\Gr^+_2(H_{e,\RR}) &
\end{tikzcd}
\end{center}
The map $\Gr^+_3(H_\RR)\to \Gr^+_2(H_{e,\RR})$ is a fiber bundle whose fiber over $[\Pi]\in \Gr^+_2(H_{e,\RR})$ is the  19-dimensional  hyperbolic  space  $\Hcal_\Pi$. The open subset  $U_e$ of  $\Gr^+_3(H_\RR)$ 
is such  that it meets each  fiber $\Hcal_\Pi$ in an open  chamber $\Kcal_\Pi$. The inequalities defining 
$U_e$ in $\Gr^+_3(H_\RR)$ are locally  finite on  $\Gr^+_3(H_\RR)$.  

We claim that the inclusion   $U_e\subset\Gr^+_3(H_\RR)$ is a 
fiberwise homotopy equivalence over $\Gr^+_2(H_{e,\RR})$.  To see this,  first construct 
a section $s$ of the projection  $U_e\to  \Gr^+_2(H_{e,\RR})$ by means of  a partition of unity.  Consider the smooth vector field on $U_e$ that is tangent to the fibers of   $U_e\to \Gr^+_3(H_\RR)$, points towards the section $s$ and whose length is the $\Gr^+_3(H_\RR)$-distance  to the section $s$.  This flow can be  integrated  to a map $[0, \infty]\times U_e\to U_e$ that gives a fiberwise deformation retraction onto $s$. As the base $\Gr^+_2(H_{e, \RR})$ is contractible, $U_e$ will be as well.

Next we observe that  $\Tcal(M)_e=U_e\cap \Tcal(M)= U_e\cap \Gr^+_3(H_\RR)^\circ$. So the locus that we remove from $U_e$ to get $\Tcal(M)_e$ is everywhere of real codimension $3$. This does not affect the simple connectivity and  it thus follows that $\Tcal(M)_e$ is simply-connected.
\end{proof}

\begin{proof}[Proof of  \ref{thm:unipotent radical}]
\label{proof:unip rad}
It remains to see that the universal bundle $f: \Ucal_{\Tcal(M)}\to \Tcal(M)$,
when restricted to $\Tcal(M)_e$,  $f_e: \Ucal_{\Tcal(M)_e}\to \Tcal(M)_e$, factors over a $\PP^1$-bundle as asserted.  
This amounts to a standard argument: we may regard $e$ as the class of a line bundle 
$\Lcal$ on $\Ucal_{\Tcal(M)_e}$ that is holomorphic on each fiber of $f_e$. This defines a linear system over $\Tcal(M)_e$ of relative rank one: 
the direct image $f_*\Lcal$ defines  a smooth complex vector bundle over $\Tcal(M)_e$  of rank two. 
If we let $\Pcal_{\Tcal(M)_e}$ stand for the projectivization of this bundle (strictly speaking of its dual, but that is the same in the rank two case), then this is a smooth $\PP^1$-bundle over $\Tcal(M)_e$ which is holomorphic on each fiber of  $f_e$. The natural map
\[
\Ucal_{\Tcal(M)_e}\to \Pcal_{\Tcal(M)_e}
\] 
then gives the fiberwise elliptic fibrations over $\Tcal(M)_e$. This completes the proof of Theorem \ref{thm:unipotent radical}.
\end{proof}

\subsection{Proof of Theorem \ref{theorem:mint1}}

In this subsection we prove most of Theorem \ref{thm:Jac}.  
Let $\Gr^+_2(H_{e,\RR})^\circ:= h''{}^{-1}\Gr^+_2(H(e)_\RR)^\circ$. We begin with the following proposition.

\begin{proposition}\label{prop:mint3}
The image of $\Tcal(M)^\circ_e$  in $\Gr^+_3(H_\RR)$ (i.e.,  the part of $\Gr^+_3(H_\RR)$ that makes $e$  the class of an integral genus one fibration) has the structure of an iterated  $C^\infty$-fiber bundle over $\Gr^+_2(H(e)_\RR)^\circ$ having $\Gr^+_2(H_{e,\RR})^\circ$  as its intermediate space; the bundle 
$\Tcal(M)^\circ_e\to \Gr^+_2(H_{e,\RR})^\circ$ is fibered into hyperbolic chambers of dimension of 19
and  $\Gr^+_2(H_{e,\RR})^\circ\to \Gr^+_2(H(e)_\RR)^\circ$ has the structure 
of a holomorphic affine line bundle.

If we divide out by the action of $H(e)$ (acting as a group of Eichler transformations \eqref{eqn:eichler}), then this yields a bundle
\[
H(e)\bs \Tcal(M)^\circ_e\to \Gr^+_2(H(e)_\RR)^\circ
\]
whose fibers have canonically the homotopy type of the $20$-dimensional  torus $T(e)=\RR/\ZZ\otimes H(e)$.
\end{proposition} 
\begin{proof} 
The first statement mostly sums up the discussion  Subsection \ref{subsection:uniprad}. 
If  we fix $[\Pi']\in \Gr^+_2(H(e)_\RR)$, then $h''^{-1}[\Pi']$ is the complex  affine line which parametrizes the 
$[\Pi]\in \Gr^+_2(H_{e,\RR})$ for which $\Pi$ maps onto $\Pi'$ and for any such $\Pi$ we have that $h'^{-1}[\Pi]$ is the 19-dimensional hyperbolic space  $\Hcal_\Pi$. In particular, $h^{-1}[\Pi]$ is contractible. 
 The action of $H(e)$ on $h^{-1}[\Pi']$ is free and properly discontinuous  and hence its  orbit space  $H(e)\bs h^{-1}[\Pi']$ is an Eilenberg-Mac Lane space for $H(e)$; it therefore has the homotopy type of the torus $T(e)$. Since $H(e)$ is abelian, there is in fact a canonical class of  homotopy equivalences between  this torus and 
$H(e)\bs h^{-1}[\Pi']$. To be precise, the projection  $H(e)\bs \Gr_3^+(H_\RR)\to  Gr^+_2(H(e)_\RR)$ induced by $h$  is fiberwise homotopy equivalent to a $T(e)$-bundle.

We observed  that $\Tcal(M)^\circ_e\to Gr^+_2(H(e)_\RR^\circ$ is the restriction of the bundle $U_e\to \Gr^+_2(H(e)_\RR$ to 
$\Gr^+_2(H(e)_\RR^\circ$ and that  its  fibers  are themselves fibered over a   complex line  with fibers hyperbolic  chambers. So in the diagram
\[
H(e)\bs \Tcal(M)^\circ_e \subset H(e)\bs h^{-1}\Gr^+_2(H(e)_\RR)^\circ\xrightarrow{h} \Gr^+_2(H(e)_\RR)^\circ
\]
the first map a fiberwise homotopy equivalence. Hence the restriction of $h$ to $H(e)\bs \Tcal(M)^\circ_e$  is fiberwise homotopy equivalent to a $T(e)$-bundle. The second assertion follows.
\end{proof}

\begin{remark}\label{rem:}
We can identify $h^{-1}[\Pi']$ with the set $K$ of $\kappa\in H_\RR$  with  $\kappa\cdot e>0$ and $\kappa\cdot\kappa=2$ as follows.  If $\Pi'$ is  a positive $2$-plane in $H(e)_\RR$, then for every $\kappa\in K$, there is unique lift of $\Pi'$ to a $2$-plane $\Pi_\kappa\subset H_{e,\RR}$ that is perpendicular to $\kappa$ (see the proof of Lemma \ref{lemma:factorization1}) 
so that then $h([\Pi_\kappa+\RR\kappa])=[\Pi']$ and $\kappa\in \Hcal_{\Pi_\kappa}$. This identifies  $h^{-1}[\Pi']$ with $K$  and does so $H(e)_\RR$-equivariantly. One may also verify directly that $K$ is contractible.
\end{remark}

\begin{corollary}\label{cor:mint}
The restriction of the commutative triangle of Theorem \ref{thm:unipotent radical} to $\Tcal(M)^\circ_e$ yields the universal bundle of marked 
integral genus one fibrations  of a K\"ahler-Einstein K3-surface:
\begin{center}
\begin{tikzcd}[column sep=scriptsize]
\Ucal_{\Tcal(M)^\circ_e }\arrow[dr] \arrow[rr]{} & &\Pcal_{\Tcal(M)^\circ_e} \arrow[dl]\\
&\Tcal(M)^\circ_e
\end{tikzcd}
\end{center}
It comes with a $\G_e$-action. Taking the quotient by that action yields the moduli space $\Mcal_\int$ of 
integral genus one fibrations of a K\"ahler-Einstein K3-surface.  

The map $\Tcal(M)^\circ_e\to \Gr_2^+(H(e)_\RR)^\circ$ induces map  of orbifolds
\begin{equation}\label{eqn:subquot}
\Mcal_\int\cong \G_e\bs \Tcal(M)^\circ_e\to \G(e)\bs \Gr_2^+(H(e)_\RR)^\circ\cong \Mcal_\int^\el,
\end{equation}
whose (orbifold) fibers have the homotopy type  of a $20$-dimensional torus and for which $\Mcal_\int^\el$ defines a section. In particular,
the  orbifold fundamental group of $\Mcal_\int$ is a split extension of the orbifold fundamental group of $\Mcal^\el_\int$ by $H(e)$. 
\end{corollary}
\begin{proof}
All these assertions follow from Proposition \ref{prop:mint3} except for the identification of $\Mcal_\int^\el$ as a section. For this recall that the K\"ahler class we assign an elliptically fibered structure on $M$ with only integral fibers and fiber class $e$ and  section class $\sigma$ is 
$\sqrt{1/2}(\sigma+3e)$. Since $(\sigma+3e)\cdot e=1$ the locus $\Mcal_\int^\el$ defines a section. 
\end{proof}

This proves Theorem \ref{theorem:mint1}.

\section{The universal Jacobian fibration}
\label{section:jacobian}

Let $X$ be a  K3-surface endowed  with holomorphic genus $1$ fibration $\pi:X\to S$ with discriminant $D$.  Recall from the introduction the  associated Jacobian fibration 
$\jac(\pi) \to S$.  Concretely, the inclusion $\ZZ_X\subset \Ocal_X$ induces a map $R^1\pi_*\ZZ\to R^1\pi_*\Ocal_X$. The Jacobian fibration is over $S^\circ:= S\ssm \supp(D)$ given as  the cokernel of this map (with the section given as the zero section) and this is all we need to know for what follows.
We shall assume that the fiber class $f\in H^2(X)$ of $\pi$ is primitive. This makes  $f^\perp/\ZZ f$ free abelian and the intersection pairing induces  a pairing  on $f^\perp/\ZZ f$. It also inherits a Hodge structure  from the one on $H^2(X)$ (as a subquotient) that is polarized by this pairing.
Let us refer to this polarized Hodge structure as the \emph{Leray-primitive subquotient of  $H^2(X)$}.

Note that if we are given a section of $\pi$  with class $\sigma\in H^2(X)$, then the Leray-primitive subquotient of  $H^2(X)$
is naturally realized as a direct summand of $H^2(X)$, namely as the orthogonal complement 
 $\ZZ\sigma+\ZZ f$. As $\sigma$ and $f$  are of type $(1,1)$,  the Hodge structure on $H^2(X)$ is then completely encoded by the one on its Leray-primitive subquotient.

For our purpose, the key result we need is the following.

\begin{proposition}\label{prop:subquotient}
Assume that all the fibers of $\pi$ are irreducible.
Then the  Leray-primitive subquotients of  $H^2(\jac(X/S))$ and $H^2(X)$ are  canonically isomorphic
as polarized Hodge structures.
\end{proposition}

The importance of Proposition \ref{prop:subquotient} is that it implies the following. 

\begin{theorem}\label{thm:jacmoduli}
The natural map 
\[\Mcal_\int\cong \G_e\bs \Tcal(M)^\circ_e\to \G(e)\bs \Gr_2^+(H(e)_\RR)^\circ\cong \Mcal_\int^\el\] in Theorem \ref{eqn:subquot} has the modular interpretation of passing to the Jacobian fibration.
This is equivalent to passing to the Leray-primitive subquotient.
\end{theorem}
\begin{proof}
Indeed, if $\pi: M\to \PP^1$ is a genus one fibration with fiber class $e$, then the polarized Hodge structure on 
$H^2(\jac(\pi))$ of the Jacobian fibration splits as the span of the fiber class $e$ and the section class (both are of type $(1,1)$) and its perp. That perp is naturally identified with Leray-primitive subquotient $e^\perp/e$. The theorem follows. \end{proof}

Theorem \ref{thm:jacmoduli} immediately implies Theorem \ref{thm:Jac}.  Proposition \ref{prop:subquotient} will follow if we succeed  in describing the Leray-primitive subquotient of  $H^2(X)$ entirely in terms of the map $R^1\pi_*\ZZ\to R^1\pi_*\Ocal_X$. This is what we will do. It will be based on 
the Leray spectral sequence for $\pi$,
\[
E^{p,q}_2=H^p(S, R^q\pi_*\ZZ)\Rightarrow H^{p+q}(X),
\]
Proposition \ref{prop:subquotient} then follows from Proposition \ref{prop:leray} below.

\begin{proposition}\label{prop:leray}
The above Leray spectral sequence degenerates on the second page.  Moreover,
\begin{enumerate}
\item  the Leray filtration on $H^2(X)$ is given by
\[
0\subset \ZZ f\subset f^\perp\subset H^2(X)
\]
with successive nonzero quotients $H^2(S)\cong\ZZ f$ and $H^1(S, R^1\pi_*\ZZ)\cong f^\perp/\ZZ f$ and 
$H^2(S,\pi_*\ZZ)\cong \ZZ$ and $H^0(S,R^2\pi_*\Zb)\cong H^2(X)/f^\perp\cong \ZZ$.

\item the  pairing on the Leray-primitive subquotient  $f^\perp/\ZZ f$ is via its identification with $H^1(S, R^1\pi_*\ZZ)$ given by the cup product 
\[
H^1(S,R^1\pi_*\ZZ)\otimes  H^1(S,R^1\pi_*\ZZ)\to H^2(S,R^2\pi_*\ZZ)\cong \ZZ,
\] 
\item
if we regard  $f^\perp/\ZZ f$ as a subquotient of $H^2(X)$ in the category of Hodge structures, then the $F^1$ of its Hodge filtration is under this isomorphism equal to the kernel of  the natural map
\[
H^1(S,R^1\pi_*\CC)\to H^1(S,R^1\pi_*\Ocal_X).
\]
\end{enumerate}
 \end{proposition}

\begin{proof}
The only possible nonzero differentials of the Leray sequence of $\pi$  are
\begin{gather*}
d_2^{0,1}: H^0(S, R^1\pi_*\ZZ)\to H^2(S,  \pi_*\ZZ)=H^2(S)\cong \ZZ \text{  and}\\
d_2^{0,2}: H^0(S, R^2\pi_*\ZZ)\to H^2(S,  R^1\pi_*\ZZ).
\end{gather*}
The kernel of $d_2^{0,1}$ contributes to $H^1(X)$ and the cokernel of $d_2^{0,2}$ contributes to
$H^3(X)$. Since $X$ is simply-connected, both cohomology groups are zero and hence these differentials are zero as well. So the Leray sequence degenerates.

The natural map $\pi^*:H^2(S)=H^2(S,  \pi_*\ZZ)\to H^2(X)$ takes the generator of  $H^2(S)$ to $f$  and so
$\ZZ f$ appears in the Leray filtration on $H^2(X)$. Since the fibers of $\pi$ are all irreducible, $R^2\pi_*\ZZ$ is the constant sheaf $\ZZ$ on $S$ and so 
$H^0(S, R^2\pi_*\ZZ)=\ZZ$. The natural map $H^2(X)\to H^0(S, R^2\pi_*\ZZ)\cong \ZZ$ is then given by integration over a fiber. Via Poincar\'e duality this is given by $a\in H^2(X)\mapsto  a\cdot f$ and hence nonzero. This also shows that $f^\perp$ is a member of the Leray filtration on $H^2(X)$.
It then follows that $H^1(S, R^1\pi_*\ZZ)\cong f^\perp/\ZZ f$. This proves (1)

Assertion (2) is a general compatibility property of the Leray spectral sequence.

For (3) we note that $R^1\pi_*\ZZ|_{S^\circ}$  is  a polarized variation of Hodge structure of weight 1, with the nontrivial member of the Hodge filtration being given by the kernel of the natural map 
$\Ocal_S\otimes R^1\pi_*\ZZ\to R^1\pi_*\Ocal_X$
restricted to $S^\circ$. The theory of polarized variation of Hodge structures (see Cox-Zucker, \cite{CZ} \S 3) affirms that then $H^1(S, R^1\pi_*\ZZ_X)$ comes with a polarized Hodge structure of weight 
$2$, which coincides the one that we get from  $f^\perp/\ZZ f$ via its identification with $H^1(S, R^1\pi_*\ZZ_X)$.
\end{proof}

\subsection{Proof of Theorem \ref{thm:nodal}}\label{subsect:proof of ref{thm:nodal}}
We first confine ourselves to  integral elliptic fibrations; the formation of the Jacobian will then enable us to lift our findings to the context of genus one fibrations.

An integral elliptic fibration admits a Weierstra\ss\ form. To be precise, fix  a $\PP^1$-bundle over $\PP^1$, denoted 
$F\to \PP^1$, endowed with a section $\sigma_0$ with self-intersection $-4$ 
(this makes it a Hirzebruch surface; the given data determine its isomorphism class).  
The  curves on $F$  that have intersection number $3$ with a fiber and $0$ 
with $\sigma_0$ make up a linear system, so are parametrized by a projective space, and its smooth members 
define a Zariski-open subset that we shall denote by $\Bscr$.

Let $B\in \Bscr$. Then $B$ and $\sigma_0$ will be  disjoint and we can form the double cover  
$X_B\to F$ ramified over the union of $\sigma_0$ and $B$. This double cover comes endowed 
with maps $\pi: X_B\to F\to \PP^1$ and a natural lift $\sigma: \PP^1\to X_B$ of $\sigma_0$ as 
well as with a (Galois) involution.  This is in fact an integral elliptic fibration (the involution acts as minus 
the identity in each fiber).  The fiber $X_{B,s}$,  $s\in \PP^1$ is smooth, nodal or cuspidal according to whether $B$ meets $F_s$ in $3$, $2$ or $1$ points. 

Conversely all integral elliptic fibrations so arise up to isomorphism. This was Miranda's point of departure for the construction of his coarse moduli spaces of elliptic fibrations \cite{miranda}:  The group $\aut(F)$ is an affine algebraic  group which it automatically preserves the fibration and the section.  It acts properly  on $\Bscr$ (so with finite isotropy groups) and if we divide  out $\Bscr$ by that action we recover $\Mcal^\el_\int$.

The condition that $X_{B,s}$ is a cuspidal fiber is equivalent  to $F_s$ meeting $B$ in a point 
with multiplicity 3. This defines in $\Bscr$ a closed, irreducible hypersurface $\Bscr_c\subset \Bscr$ whose generic point parametrizes the $X_B$' s that have precisely one such fiber. Let 
$j:\Delta \to \Bscr$ be a holomorphic map from the complex unit disk  such that $j^*\Bscr_c=(0)$ as divisors. This  defines a family  
$\pi:\Xcal\to \Delta\times\PP^1$ of integral elliptic fibrations. All singular fibers of $\Xcal\to \Delta\times\PP^1$ 
are nodal, save for one over $(0,s)$ for some $s\in \PP^1$. If $o\in X_{0,s}$ is the cusp, then we find that after a 
suitable coordinate change on $\Delta$ and a choice of a suitable local coordinate $u$ of $\PP^1$ at $s$, this family takes near  $(0,s)$ the simple Weierstra\ss\ form
\[
y^2=x^3- 3tx +2u
\]

The  critical points of the map $\pi$ are given by 
$(x,y,t,u)=(\pm \sqrt{t},0,t,\pm t\sqrt{t})$ on which $\pi$ takes the value 
$\mp 2t\sqrt{t}+2u$. In other words, for $t\not=0$, we find two nodal fibers over the point satisfying  
$u^2=t^3$. If we let $t$ traverse counterclockwise a circle $|t|=\eps$ with $\eps>0$ small, then the 
two critical values trace out the third power of a simple braid. This settles the analogue of Theorem 
\ref{thm:nodal} in $\Mcal^\el_\int$.

Theorem \ref{thm:nodal} itself (i.e., the corresponding result for $\Mcal_\int$),  then follows via Corollary \ref{cor:mint}.

\section{Applications to the mapping class group}
\label{section:applics}
Our results on moduli spaces have consequences for the  mapping class group of $M$. In order to state these, let us first agree on the following terminology and notation. Noting that the group $\Diff (M)\times\Diff (\PP^1)$ acts on the set of (smooth) genus one fibrations 
$\pi:M\to \PP^1$ that  have  24 nodal fibers
\begin{equation}
\label{eq:action8}
(h,h')\cdot \pi:= h'\circ \pi\circ  h^{-1},
\end{equation}
then we say that  $(h,h')$ is an \emph{automorphism 
of $\pi$}  if $(h,h')$ fixes $\pi$.  Denote the stabilizer  of $\pi$ in $\Diff (M)\times\Diff (\PP^1)$ by $\Diff(\pi)$. If in addition $h'$ is the 
identity (so that $h$ is simply a diffeomorphism of $M$ that preserves each fiber of $\pi$), then we 
say that $h$ is an \emph{automorphism over $\PP^1$}.  We denote the subgroup of such 
elements by  $\Diff (M /\PP^1)$.  We denote the corresponding component groups accordingly by $\Mod(\pi)$ resp.\  $\Mod(M/\PP^1)$.  

So $\Mod(M/\PP^1)$ is the kernel  of the evident  forgetful homomorphism from $\Mod(\pi)$ to the spherical braid group on 24 strands. Such mapping classes may arise as monodromies: if  $\pi:M\to \PP^1$ shows up as a fiber (over $o\in T$, say) in a family of genus one fibered K3 manifolds
\begin{center}
\begin{tikzcd}[column sep=scriptsize]
\Ucal\arrow[dr] \arrow[rr]{} & &\Pcal \arrow[dl]\\
&T
\end{tikzcd}
\end{center}
where $\Pcal\to T$ is a $\PP^1$-bundle, and the relative discriminant $\Dcal$ of $\Ucal\to \Pcal$  is an unramified degree 24 cover of $T$, then such a family is locally trivial in the smooth category. 
So its structural group is $\Diff (M)\times\Diff (\PP^1)$ and the monodromy of this family is a homomorphism $\pi_1(T,o)\to \Mod(\pi)$. If the degree $24$ cover is trivial, then a trivialization of that cover extends to $\Pcal$ and hence we can take the structural group to be $\Diff (M)$. The monodromy group then takes its values in  
$\Mod(M/\PP^1)$.

\begin{corollary}\label{cor:Gamma_epsilon_realization}
The following hold:
\begin{enumerate}
\item The set of holomorphic genus one fibrations $\pi:M\to\Pb^1$ with primitive fiber class and 
$24$ singular nodal fibers, lie in a single  orbit of $\Diff (M)\times\Diff (\PP^1)$ relative the action given by \eqref{eq:action8}.

\item If $\pi:M\to \PP^1$ is such  a genus one fibration with fiber class $e$, then the image of $\Diff(\pi)$ in $\Gamma$ under the map $\Diff(M)\to \Gamma\subset\Orth(H)^+$ is all of $\Gamma_e$.
\end{enumerate}

\end{corollary}

\begin{proof}
Moishezon proved that, up to the $\Diff (M)\times\Diff (\PP^1)$ 
action, any such $\pi:M\to\Pb^1$ can be moved into a specific elliptic fibration, with singular fibers at $24$th roots of unity, and monodromies alternating in clockwise order between the two standard unipotent generators of $\SL(2,\Zb)$; see \cite{FM}, Theorem 3.6, page 169. This implies the first claim.

To prove the second claim, note that there is 
a universal family over  $\Mcal_\nod$, in the sense of  our conventions on moduli spaces in \S\ref{conv:moduli}.  The fact that $\Mcal_\nod$ is connected as a moduli space implies the claim.
\end{proof}

\subsection{Fiber-preserving diffeomorphisms}
Recall that an \emph{affine structure} on manifold is specified by an atlas whose coordinate changes  are affine-linear and which is maximal for that property; this is equivalent to give  a flat, torsion free  connection on its   tangent bundle.

We begin with the observation that if  $E$ is a closed genus $1$ surface  endowed with an affine structure, then the identity component   of the automorphism group of $E$ is isomorphic to the  $2$-torus  $T^2\cong (\RR/\ZZ)^2$ for which $E$ is a torsor. Denote that group by $\Tr (E)$ and refer to it as the \emph{translation group} of $E$. A  complex structure on $E$ determines such an affine structure because there is then a unique flat metric compatible with the complex structure that gives that fiber unit volume.

The situation is similar for a once or twice-punctured 2-sphere: an affine structure 
makes such a surface a torsor of  the identity component  of its automorphism 
group, which is then isomorphic  to the additive group resp.\ the multiplicative group of 
the complex field.  

Let $\pi :M\to \PP^1$ be a genus one fibration with only nodal fibers as singular fibers and with discriminant $D$. The complex structure determines in  each fiber an affine structure (depending smoothly on the base point). So the structure group of $\pi$, at least over $\PP^1\ssm D$,  is the 
semi-direct product  $T^2\rtimes \SL_2(\ZZ)$. This determines a subgroup 
$\Diff_\aff (\pi)\subset \Diff(\pi)$ as the subgroup of fibration-preserving diffeomorphisms that also preserve this fiberwise-affine structure. We put 
\begin{gather*}
\Diff_\aff(M/\PP^1):=\Diff(M/\PP^1) \cap \Diff_\aff (\pi),\\ 
\Tr(M/\PP^1):=\{\text{fiberwise translations}\}\\
\end{gather*}
Note that $\Tr(M/\PP^1)$ is an abelian subgroup of $\Diff^+_\aff(M/\PP^1)$ that is normal in $\Diff^+_\aff (\pi)$.
We call its connected component group the \emph{smooth Mordell-Weil} group of $\pi$, denoting it
\[
 \MW(M/\PP^1):=\pi_0(\Tr(M/\PP^1)).
 \] 
 
We could also include fiberwise involutions (acting as minus the identity in a 2-torus).  
 Since these exist globally,  this gives us a semi-direct product $\Tr(M/\PP^1)\rtimes\mu_2$ contained in
$\Diff_\aff(M/\PP^1)$. We have a corresponding semi-direct product  $\MW(M/\PP^1)\rtimes\mu_2$ in $\Diff_\aff(M/\PP^1)$.

For any of the other diffeomorphism groups, its   connected component group  
will be denoted  in the usual manner (replace $\Diff$ by $\Mod$). 

If $\pi : M\to \PP^1$ appears as a member of  a family $\Ucal\to \Pcal\to T$ of such fibrations whose discriminant is locally constant, then the monodromy defines an element of $\Mod_\aff(\pi)=\pi_0(\Diff_\aff(\pi))$. 

\begin{theorem}\label{thm:mapping classes preserving a fibration}\label{thm:nielsenMW}
Let $\pi :M\to \PP^1$ be a  genus one fibration with primitive  fiber class $e$ and $24$ singular, nodal fibers. Then the action of $\Mod_\aff(\pi)$ on $H=H^2(M)$ has image $\G_e$. 

Moreover, the smooth Mordell-Weil group  $\MW(M/\PP^1)$ maps onto the unipotent radical $R_u(\G_e )$ of $\G_e$ (which we recall, can be identified with the rank 20 lattice  $H(e)=e^\perp/\ZZ e$) and there is a Nielsen realization for this subgroup in the sense that the map
$\Diff^+_\aff(M/\PP^1)\to \MW(M/\PP^1)$ admits  a section homomorphism. This section extends to
the semidirect product with $\mu_2$ giving a group homomorphism $R_u(\G_e )\rtimes \mu_2\to \Diff(M/\PP^1)$.
\end{theorem}
\begin{proof}
The first part of the theorem follows from  Theorem \ref{thm:unipotent radical} and the second part from Corollary \ref{cor:MWNielsen}; see also the discussion after this corollary.
\end{proof}

\bigskip{\noindent
Dept. of Mathematics, University of Chicago\\
bensonfarb@gmail.com\\ 
e.j.n.looijenga@uu.nl

\end{document}